%% file: Z2alt-AG-ostp_arX4.tex
\documentclass[11pt]{article}

\title{\vskip-1.5em\sc Constructing alternating $2$-cocycles on Fourier algebras}
\author{\sc Y. Choi}
\date{17th April 2021}

\usepackage[a4paper,left=27mm,right=27mm,top=30mm,bottom=33mm]{geometry}
\usepackage{sectsty}
\allsectionsfont{\sffamily}
\usepackage[colorlinks=true,linkcolor=blue,citecolor=red]{hyperref}

\usepackage[dvipsnames]{xcolor}

\input{zaltagostp_mac}

\input{widehack}


\renewcommand{\dt}[1]{\textsf{\textcolor{Bittersweet}{#1}}}

\begin{document}

\maketitle

\begin{abstract}
Building on recent progress in constructing derivations on Fourier algebras, we provide the first examples of locally compact groups whose Fourier algebras support non-zero, alternating $2$-cocycles; this is the first step in a larger project. Although such $2$-cocycles can never be completely bounded, the operator space structure on the Fourier algebra plays a crucial role in our construction, as does the opposite operator space structure.

Our construction has two main technical ingredients: we observe that certain estimates from \cite{LLSS_WAAG} yield derivations that are ``co-completely bounded'' as maps from various Fourier algebras to their duals; and we establish a twisted inclusion result for certain operator space tensor products, which may be of independent interest.


(primary);
43A30,
46M05
(secondary)

\end{abstract}


\begin{section}{Introduction}\label{s:intro-motivation}

\begin{subsection}{Background context and our main application}
Fourier algebras of locally compact groups have been a fertile source of examples in the study of general Banach function algebras, while also having some important applications to the study of operator algebras associated to group representations (see e.g.\ \cite{HK_AP-group}).
One theme with a long history is the study of how properties of a group $G$ are reflected in properties of its Fourier algebra $\FA(G)$.
For instance, if $G$ is compact and non-abelian and $f\in\FA(G)$, the matrix-valued Fourier coefficients $\widehat{f}(\pi)$ must decay at a certain rate as $\pi$ ``tends to infinity'', which intuitively suggests that $f$ should have a degree of differentiability or H\"older continuity.
This heuristic underlies a theorem of Johnson\footnotemark\
that when $G$ is $\SO(3)$ or $\SU(2)$, there are non-zero derivations from $\FA(G)$ to its dual; using general restriction theorems for Fourier algebras, it follows that the same is true for any $G$ that contains a closed copy of $\SO(3)$ or $\SU(2)$.
\footnotetext{For details, see the proof of \cite[Theorem~7.4]{BEJ_AG} or the discussion in \cite[\S3]{CG_WAAG1}.}

Johnson's result went counter to some expectations at the time: for given a Lie group $G$ and some $X$ in its Lie algebra, the Lie derivative along $X$, viewed as a continuous operator on $C_c^\infty(G)$, does not extend to a continuous map $\FA(G)\to C_b(G)$. (If such an extension existed it would yield non-zero continuous point derivations on $\FA(G)$, contradicting the fact that points of $G$ are sets of synthesis for~$\FA(G)$.)
Enlarging the codomain from $C_b(G)$ to $\FA(G)^*$ allows $\partial_X f$ to be a distribution rather than a function, but a separate averaging argument is needed to explain why $f\mapsto \partial_X f$ has any chance of being continuous from $\FA(G)$ to $\FA(G)^*$.

Despite these technical difficulties, 
Johnson's result has been greatly extended in recent years, by using the operator-valued Fourier transform for certain Type I groups to make explicit calculations and estimates:
 see the papers \cite{CG_WAAG1}, \cite{CG_WAAG2} and \cite{LLSS_WAAG}.
These papers were motivated by a conjecture posed by Forrest and Runde in \cite{ForRun_amenAG}, which predicted exactly which groups $G$ allow non-zero derivations $\FA(G)\to\FA(G)^*$; the conjecture was confirmed for all Lie groups in~\cite{LLSS_WAAG}, and a recent preprint of Losert~\cite{Los_WAAG_pre} contains a solution in full generality.

For any commutative Banach algebra $A$, the space of derivations $A\to A^*$ coincides with the first Hochschild cohomology group $\cH^1(A,A^*)$.
The higher-degree groups $\cH^n(A,A^*)$ potentially capture more information about~$A$, but have proved to be extremely difficult to calculate except in degenerate cases.
A more promising approach to finding computable invariants,
which was first pursued systematically in \cite{BEJ_ndimWA},
is to consider the \emph{alternating part} of $\cH^n(A,A^*)$; alternating cocycles are more tractable than general ones since they are built from derivations (and there is also conceptual motivation for singling out this class, see Remark~\ref{r:HKR behind the scenes} below).
Nevertheless, the existence of non-zero alternating cocycles on a Banach function algebra is very sensitive to properties of the given norm, and is not guaranteed by simply having ``enough derivations'', as illustrated in Example~\ref{eg:lip-alg} below.

Given the recent progress in studying derivations on Fourier algebras, it is natural to turn our attention to alternating cocycles on~$\FA(G)$.
The present paper is the first step in getting this larger programme off the ground, by producing the first examples of groups whose Fourier algebras support non-zero alternating $2$-cocycles.
In fact, we show that not only do such groups exist, but they occur in abundance among the classical Lie groups.

\begin{thm}
\label{t:headline}
Let $n\geq 4$ and let $G$ be one of the groups $\SU(n)$, $\SL(n,\Real)$ or ${\rm Isom}(\Real^n)$. Then there is a non-zero, continuous, alternating $2$-cocycle on $\FA(G)$.
\end{thm}

Theorem~\ref{t:headline} is a special case of a much more general result, which in turn follows by combining Theorems \ref{t:mainthm} and~\ref{t:underyournoses} below.
The starting point for the proof of Theorem~\ref{t:mainthm} is the canonical identification of $\FA(H\times L)$ with the operator space projective tensor product of $\FA(H)$ and $\FA(L)$, valid for any locally compact groups $H$ and~$L$. However, it should be emphasised that our proof requires more than a merely formal use of operator spaces and completely bounded maps; one can show that there are no completely bounded, non-zero, alternating cocycles on Fourier algebras.

The key extra ingredient in the proof of Theorem~\ref{t:mainthm} is the following surprising phenomenon, which may have independent interest for those working in operator space theory.

\begin{thm}[``Twisted inclusion'']
\label{t:twisted inclusion}
Let $X$ and $Y$ be operator spaces, and let $\til{Y}$ denote $Y$ equipped with its \dt{opposite operator-space structure}. Let $\ptp$ and $\itp$ denote the projective and injective tensor products of operator spaces. Then
\[ \norm{w}_{X\smallitp\til{Y}} \leq \norm{w}_{X\ptp Y} \qquad\text{for all $w\in X\tp Y$.} \]
That is: the identity map on $X\tp Y$ extends to a contractive linear map $\theta_{X,Y}: X\ptp Y\to X\itp \til{Y}$.
\end{thm}

The opposite operator space structure may be thought of, intuitively, as a ``mirror image'' of the given one.
One of the background aims of this paper is to argue that in applying operator-space methods to the study of Fourier algebras, it may be useful to work simultaneously with both the canonical operator space structure on $\FA(G)$ and its mirror image.

\end{subsection}

\begin{subsection}{An overview of the main technical difficulties}
\label{ss:intro-justify}
We now sketch why the natural attempt to prove Theorem~\ref{t:headline}, by merely following algebraic recipes in an appropriate functional-analytic category, does not work, and indicate the new ideas needed to overcome these difficulties. None of the material here is logically necessary for the proof of Theorem~\ref{t:headline} or Theorem~\ref{t:twisted inclusion}, and it may be skipped if the reader wishes to get straight to the precise mathematical details.

We start with a simple example from commutative algebra.

\begin{eg}\label{eg:ur}
Consider the $\Cplx$-algebra $\Cplx[z,w]$ and define $F_0:\Cplx[z,w]\times\Cplx[z,w]\to\Cplx[z,w]$ by
\[
F_0(f_1,f_2) =
 \frac{\partial f_1}{\partial w}\frac{\partial f_2}{\partial z} - \frac{\partial f_2}{\partial w}\frac{\partial f_1}{\partial z} \qquad(f_1,f_2\in \Cplx[z,w]). \]
 Then $F_0$ is an alternating $2$-cocycle (on $\Cplx[z,w]$ with coefficients in $\Cplx[z,w]$). Note that $\Cplx[z,w]\cong\Cplx[z]\tp \Cplx[w]$, and we may view $\partial/\partial z$ as $(d/dz) \tp \iota$ and $\partial/\partial w$ as $\iota\tp(d/dw)$.
 \end{eg}
 
This is a special case of a general algebraic construction: given two commutative $\Cplx$-algebras $\sA$ and $\sB$, a symmetric $\sA$-bimodule~$\sX$ and a symmetric $\sB$-bimodule~$\sY$, and derivations $D_{\sA} : \sA \to \sX$ and $D_{\sB}: \sB\to \sY$, we can ``wedge together'' the two amplified derivations
\[
D_{\sA}\tp\iota_{\sB}: \sA\tp \sB\to\sX\tp \sB
 \quad\text{and}\quad
 \iota_{\sA}\tp D_{\sB}: \sA\tp\sB \to \sA\tp\sY
\]
to obtain an alternating $2$-cocycle $F_0 : (\sA\tp \sB)\times(\sA\tp \sB) \to \sX\tp \sY$. Under mild conditions on $D_{\sA}$ and $D_{\sB}$, $F_0$ will be non-zero.

Therefore, in cases where we have non-zero derivations from $\FA(H)$ and $\FA(L)$ into appropriate modules, one could try to obtain alternating $2$-cocycles on $\FA(H\times L)$ by applying the natural analogue of this construction for the category of Banach spaces. Unfortunately this would only yield a non-zero $2$-cocycle on the Banach space projective tensor product $\FA(H)\bsptp\FA(L)$, and the natural map $\FA(H)\bsptp\FA(L)\to \FA(H\times L)$ is never surjective in such cases.
On the other hand, $\FA(H\times L)$ can be identified with the operator-space projective tensor product $\FA(H)\ptp\FA(L)$. But in the natural analogue of the algebraic construction for the category of operator spaces, one must start with completely bounded derivations from $\FA(H)$ and $\FA(L)$ into symmetric ``c.b.-bimodules'', and it was shown independently by Spronk and Samei that the only such derivations are identically zero: see \cite[Theorem 5.2]{Sam_cb-loc-der} or~\cite{Spr_OWA}.

To resolve this apparent stalemate we need a new approach, which is to examine what occurs if we start from a pair of non-zero derivations $D_H:\FA(H)\to\FA(H)^*$ and $D_L:\FA(L)\to\FA(L)^*$ that become completely bounded \emph{after changing the operator space structure} on the codomains.
Although this breaks certain aspects of the algebraic construction, something survives if $D_H$ and $D_L$ are completely bounded when $\FA(H)^*$ and $\FA(L)^*$ are equipped with their opposite operator space structures (it is an unrecorded observation of the author and Ghandehari that this property holds for the derivations constructed in~\cite{CG_WAAG1}).
For then, as we shall show in Proposition~\ref{p:twist-wedge}, combining $D_H$ and $D_L$ and doing some careful book-keeping yields
 a bounded (but not completely bounded) bilinear map
\begin{equation}\label{eq:almost there}
F: \FA(H\times L) \times \FA(H\times L) \longrightarrow \FA(H)^*\ptp\longopp{[\FA(L)^*]}
\end{equation}
which behaves like an alternating $2$-cocycle on the dense subalgebra $\FA(H)\tp\FA(L)$.

This is still not enough to obtain Theorem~\ref{t:headline}, since the right-hand side of \eqref{eq:almost there} is not an $\FA(H\times  L)$-bimodule in the cases where $D_H$ and $D_L$ exist with the required properties. The saving grace is Theorem~\ref{t:twisted inclusion}, which allows us to embed this space continuously (but not completely boundedly!) into $\FA(H\times L)^*$.
Note that in general $X\ptp Y$ and $X\itp\til{Y}$ are \emph{incomparable as operator spaces} (just take $X=\Cplx$),
and so Theorem \ref{t:twisted inclusion} does not seem to follow just from the universal/extremal properties of $\ptp$ and $\itp$ in the category of operator spaces.
Instead, our proof proceeds by embedding $X$ and $Y$ into $B(E)$ and $B(F)$ for Hilbert spaces $E$ and $F$, which allows us to calculate or bound various tensor norms on $B(E)\tp B(F)$ by viewing elements of this space as elementary operators on Schatten classes.

It is striking that to prove Theorem~\ref{t:headline}, which on the face of it makes no reference to operator spaces and completely bounded maps, we are driven to make substantial use of such techniques.

\end{subsection}

\begin{subsection}{Structure of the paper}
Let us now describe the organization of the rest of this paper.
In Section~\ref{s:prelim} we establish some global conventions for our notation, and set up the definition of $\FA(G)$ that is most suitable for this paper.
In Section~\ref{s:Z2alt} we give the key definitions of derivations and alternating $2$-cocycles for commutative Banach algebras, illustrating the general definitions with some key examples that will motivate our proof of Theorem~\ref{t:mainthm}. We also record some basic constructions that were not mentioned in \cite{BEJ_ndimWA}, as they may be useful for subsequent work.

Section~\ref{s:UNNAMED} has two purposes.
We introduce the key notion of a \dt{co-completely bounded map} (co-cb for short) between two given operator spaces; and we collect some results concerning the canonical operator space structure on $\FA(G)$, some of which are only stated implicitly in the literature.
 In doing so, we spend time on the crucial notion of the \dt{opposite operator space structure} and its functorial properties; this requires us to set out some basic properties that do not seem to be mentioned explicitly in \cite{ER_OSbook} or~\cite{Pis_OSbook}.

Section~\ref{s:make cocycles} contains the main work needed to establish Theorem~\ref{t:headline}.
En route, we give the proof of Theorem~\ref{t:twisted inclusion}, reducing the problem to a special case which is handled by means of an interpolation argument.
The section ends by stating and proving the main technical theorem of this paper,
Theorem~\ref{t:mainthm}, which says loosely speaking that we can construct non-zero $2$-cocycles given enough non-zero co-cb derivations.

For some groups where non-zero derivations have been constructed, the co-cb property can be read out of the explicit formulas in~\cite{CG_WAAG1}; but in fact those cases and more besides can be obtained by repurposing some technical results from \cite{LLSS_WAAG}. Details are given in Section~\ref{s:co-cb-der-AG}, culminating in Theorem~\ref{t:underyournoses} which provides co-cb derivations in all cases needed to establish Theorem~\ref{t:headline}.
Finally, in Section~\ref{s:finalremarks} we make some remarks and pose some questions, with a view to possible avenues for future work.

We have attempted to make this paper accessible to workers in the general area of Banach algebras, or functional analysts interested in structural properties of particular Banach algebras. 
In particular, we have not assumed any prior familiarity with either Fourier algebras of locally compact groups, or the Hochschild cohomology groups of Banach algebras, and have tried to include a small amount of extra motivation for these objects.
On the other hand, we do assume some previous exposure to the basic language of operator spaces and completely bounded maps.

\end{subsection}

\end{section}

\begin{section}{Preliminaries}
\label{s:prelim}

\begin{subsection}{Conventions and notation}
Throughout this article, all derivations and cocycles from Banach algebras into Banach bimodules are tacitly assumed to be norm-continuous. This is the convention adopted, for instance, in \cite{BEJ_ndimWA}, and this article will not be concerned with any issues of automatic continuity.

The algebraic tensor product of two complex vector spaces $E$ and $F$ is denoted by $E\tp F$. 
The term ``map'' is our short-hand for ``linear map'' or ``linear operator''.

All Banach spaces are defined over complex scalars.
For a Banach space $E$, $B(E)$ denotes the algebra of bounded linear operators on $E$.
The adjoint of a bounded linear map $f:E\to F$ between Banach spaces is denoted by $f^*:F^*\to E^*$, with one important exception: when $E$ and $F$ are Hilbert spaces, we shall denote the adjoint map $F^*\to E^*$ by $f^\#$ to avoid confusion with the adjoint in the sense of operators between Hilbert spaces.

One slight departure from usual conventions is that when $E$ is a Hilbert space, we shall formulate various constructions in terms of the dual space $E^*$ rather than the conjugate space~$\overline{E}$; of course the two spaces are canonically isomorphic as Banach spaces via the Riesz--Fr\'echet theorem.
(This decision is motivated by issues concerning operator space structures, but is ultimately only a matter of notational preference.)
When $E=L^2(\Om,\mu)$ for some measure space $(\Om,\mu)$ and $\eta\in L^2(\Om,\mu)$ we shall write $\evetabar$ for the functional
\begin{equation}\label{eq:L^2-iso-to-dual}
 \evetabar: \xi\mapsto \ip{\xi}{\overline{\eta}}_{L^2(\Om,\mu)}=\int_{\Om} \xi\eta \,d\mu
\end{equation}
so that $\eta\mapsto \evetabar$ is a \emph{linear} isomorphism of Banach spaces $L^2(\Om,\mu)\to L^2(\Om,\mu)^*$.

The projective tensor product of Banach spaces $E$ and $F$ is denoted by $E\bsptp F$; the Hilbertian tensor product of Hilbert spaces $V$ and $W$ is denoted by $V\tp_2 W$.

Our notational conventions for operator spaces and completely bounded maps will be set out in Section~\ref{s:UNNAMED}, since they are not needed until Section~\ref{s:make cocycles}.
\end{subsection}

\begin{subsection}{Fourier algebras of locally compact groups}
Fourier algebras originated from the study of $L^1$-group algebras of locally compact \emph{abelian} (LCA) groups. Given a LCA group $\Gamma$ with Pontrjagin dual $G=\widehat{\Gamma}$, the Fourier transform $\cF:L^1(\Gamma)\to C_0(G)$ is an injective algebra homomorphism, and so one may study the convolution algebra $L^1(\Gamma)$ by examining the function algebra $\cF(L^1(\Gamma))$ equipped with the norm pushed forwards from $L^1(\Gamma)$. This function algebra, denoted by $\FA(G)$, is now known as the Fourier algebra of~$G$.

Using Bochner's theorem, one can characterize $\FA(G)$ in terms of positive-definite functions on $G$, without reference to the group $\Gamma$. Guided by this philosophy, and following work of Godement and Stinespring in the unimodular case,
Eymard \cite{Eym_BSMF64} gave a definition of $\FA(G)$ that is valid for any locally compact group~$G$;
details of the foundational results from Eymard's original paper, and much more besides, may be found in the recent book \cite{KL_AGBGbook}.
However, for most of this paper, it is more convenient to work with an alternative description of~$\FA(G)$. The one we use is also standard, and can be found in e.g.~\cite[Defn.~VII.3.8]{Takesaki_vol2}, but our presentation has some cosmetic differences from the usual one since we wish to work with duals of Hilbert spaces rather than their conjugates.

Fix a choice of left Haar measure on $G$, denoted by $ds$, and let $\lambda: G \to \cU(L^2(G))$
denote the left regular representation of $G$ on $L^2(G)$
 defined by $[\lambda(x)f](s) = f(x^{-1}s)$.
Given $\xi,\eta\in L^2(G)$ and $x\in G$, let
\begin{equation}\label{eq:left-reg coeff}
\Psi(\xi\tp\evetabar)(x) =  \ip{\lambda(x)\xi}{\overline{\eta}}_{L^2(G)} = \int_G \xi(x^{-1}s)\eta(s)\,ds \,.
\end{equation}
This defines a contractive linear map $\Psi: L^2(G)\bsptp L^2(G)^*\to C_0(G)$ whose range we denote by $\FA(G)$.
We equip $\FA(G)$ with the quotient norm of $L^2(G)\bsptp L^2(G)^*/\ker(\Psi)$.
Using Fell's absorption theorem, it can be shown that $\FA(G)$ is closed under pointwise product and the norm on $\FA(G)$ is submultiplicative.
(See e.g.~\cite[\S4.1]{Zwarich_MSc} for a quick exposition of these results.)

 Thus $\FA(G)$ is a Banach algebra of functions on $G$, called the \dt{Fourier algebra} of $G$.
The equivalence of this definition with the original one can be extracted from the results in \cite[Chapitre~3]{Eym_BSMF64} (see also \cite[Prop.~2.3.3]{KL_AGBGbook}).
Note that with our definition, if $G$ is a LCA group then we recover the isomorphism $L^1(\widehat{G})\cong\FA(G)$ using Parseval's theorem.

Given Banach algebras $A$ and $B$ and a continuous homomorphism $A \to B$,  derivations and cocycles on $B$ can be pulled back to give derivations and cocycles on $A$ (a precise statement will be given in Lemma~\ref{l:pull back}).
It is therefore useful to identify homomorphic images of $\FA(G)$ which have a simpler form, so that we can build derivations or cocycles on those algebras instead. The following result was proved by Herz in a more general setting, and is usually known as \dt{Herz's restriction theorem}.

\begin{thm}[Herz; McMullen]
If $G$ is a locally compact group and $G_1$ is a closed subgroup, then restriction of functions $C_0(G)\to C_0(G_1)$ defines a norm-decreasing homomorphism $\FA(G)\to\FA(G_1)$ which is a quotient map of Banach spaces.
\end{thm}

For a detailed proof and some historical comments, see~\cite[\S2.6 and \S2.10]{KL_AGBGbook}.
It is worth noting that if we wish to exploit the restriction theorem to construct derivations or non-trivial cocycles, we have to take $G_1$ to be non-abelian; Fourier algebras of abelian groups are \emph{amenable} (in the Banach-algebraic sense) and hence all cohomology with dual-valued coefficients vanishes.

In concrete cases, such as those in Theorem~\ref{t:headline}, once we have constructed a well-defined cocycle on a particular Fourier algebra it will be obvious from context that this cocycle is not identically zero. However, in order to state our main technical theorem (Theorem~\ref{t:mainthm}) in its greatest generality, it will be convenient to use the following lemma as a ``soft'' work-around.

\begin{lem}\label{l:dense-powers}
Let $G$ be a locally compact group. Then $\{ a^4 \colon a\in \FA(G)\}$ and $\{b^2\colon b\in \FA(G)\}$ both have dense linear span in $\FA(G)$.
\end{lem}

\begin{proof}
Let $A_0 = \FA(G)\cap C_c(G)$. $A_0$~is a dense subalgebra of $\FA(G)$ (this is immediate from Eymard's original definition, but also follows easily from the one given above). Moreover, for each $f\in A_0$ there exists $g\in A_0$ such that $fg=f$; see \cite[Lemme~3.2]{Eym_BSMF64} or \cite[Prop.\ 2.3.2]{KL_AGBGbook}.
Hence, by the usual polarization identity
\[ ab = \frac{1}{4} \left[ (a+b)^2 - (a-b)^2\right]\;, \]
we have $A_0 \subseteq \lin\{ b^2\colon b\in A_0\}$. (The converse inclusion also holds.)
Similarly, by using the identity
\[
x^2y^2 = \frac{1}{24} \left[ (x+y)^4 + (x-y)^4 - (x+iy)^4 - (x-iy)^4 \right]\,,
\]
the same ``absorption trick'' as above yields $A_0 \subseteq \lin\{ b^2\colon b\in A_0\} \subseteq \lin \{a^4\colon a\in A_0\}$. (Once again the converse inclusion is trivial.)
Thus, both
 $\lin\{ a^4 \colon a\in \FA(G)\}$ and $\lin\{b^2\colon b\in \FA(G)\}$ contain $A_0$ and hence are dense in $\FA(G)$.
\end{proof}

\end{subsection}

\end{section}

\begin{section}{Alternating $2$-cocycles on commutative Banach algebras}\label{s:Z2alt}

\begin{subsection}{Definitions and preliminaries}
We assume familiarity with the basic language of Banach algebras and Banach bimodules over them.
Let $A$ be a Banach algebra and $X$ a Banach $A$-bimodule. For $n\geq 0$ let $\cC^n(A,X)$ be the space of bounded $n$-multilinear maps $A\times\dots\times A\to X$, with the convention that $\cC^0(A,X)=X$. There are maps $\delta^n:\cC^n(A,X)\to \cC^{n+1}(A,X)$, called the \dt{Hochschild coboundary operators}, that satisfy $\delta^{n+1}\circ\delta^n=0$ for all $n\geq 0$. We only need the cases $n=1$ and $n=2$, which have the following explicit form:
\begin{itemize}
\item[--]
for $T\in \cC^1(A,X)$ and $a,b\in A$, $\delta^1 T(a,b) \defeq a\cdot T(b)-T(ab)+T(a)\cdot b$;
\item[--]
for $T\in \cC^2(A,X)$ and $a,b,c\in A$, $\delta^2 T(a,b,c) = a\cdot T(b,c) - T(ab,c) + T(a,bc) - T(a,b)\cdot c$.
\end{itemize}
Let $\cZ^n(A,X)\defeq\ker\delta^n$;
elements of this space are called \dt{$n$-cocycles.} Note that $1$-cocycles are the same as \dt{derivations.}
Let $\cB^n(A,X)=\im\delta^{n-1}$;
elements of this space are called \dt{$n$-coboundaries.} The quotient space $\cH^n(A,X)\defeq \cZ^n(A,X)/\cB^n(A,X)$ is the \dt{$n$th Hochschild cohomology group.}

For the rest of this section, we only consider commutative Banach algebras $A$, and those Banach $A$-bimodules $X$ which are \dt{symmetric} in the sense that $a\cdot x = x\cdot a$ for all $a\in A$ and $x\in X$.
Note that if $A$ is commutative, then both $A$ itself and its dual $A^*$ are symmetric Banach $A$-bimodules.

Let $S_n$ denote the symmetric group on $\{1,\dots,n\}$, and let $\sigma \mapsto (-1)^\sigma$ denote the signature homomorphism $S_n\to \{\pm 1\}$.
A given $T\in \cC^n(A,X)$ is called \dt{symmetric} if
\[
T(a_1,\dots,a_n) = T(a_{\sigma(1)}, \dots, a_{\sigma(n)})
\qquad\text{for all $a_1,\dots, a_n\in A$ and $\sigma\in S_n$}
\]
and \dt{alternating} if
\[
T(a_1,\dots,a_n) = (-1)^\sigma T(a_{\sigma(1)}, \dots, a_{\sigma(n)})
\qquad\text{for all $a_1,\dots, a_n\in A$ and $\sigma\in S_n$.}
\]
We write $\cZ^n_{\rm alt}(A,X)$ for the space of all alternating $n$-cocycles.

\begin{rem}
In the case $n=2$, every $T\in\cC^2(A,X)$ is the sum of a symmetric part and an alternating part. Also, every $2$-coboundary is symmetric (since $A$ is commutative and $X$ is symmetric). Thus the natural map $\cZ^2_{\rm alt}(A,X)\to \cH^2(A,X)$ is actually an \emph{injection}.
\end{rem}

Following \cite[Definition 2.2]{BEJ_ndimWA}, we say that $T\in \cC^n(A,X)$ is an \dt{$n$-derivation} if it is a derivation in each variable separately. If $T$ is either symmetric or alternating, then to verify the $n$-derivation property it suffices to check that $T$ is a derivation in the first variable, i.e.\ to check the identity
\[ T(bc,a_2,\dots, a_n) = b\cdot T(c,a_2,\dots,a_n) + T(b, a_2,\dots, a_n)\cdot c \qquad\text{for all $b,c,a_2,\dots a_n\in A$.}
\]
Given our earlier definitions, a straightforward calculation shows that every $2$-derivation is a $2$-cocycle.
(It is crucial here that $A$ is commutative and $X$ is symmetric.)
In particular, every alternating $2$-derivation defines an element of $\cZ^2_{\rm alt}(A,X)$.

The alternating $2$-cocycles that we shall construct when proving Theorem~\ref{t:headline} are created as alternating $2$-derivations.
For sake of completeness, we mention that every alternating $2$-cocycle turns out to be a derivation in the first variable, and hence (by the remarks above) is an alternating $2$-derivation. We omit the details, since this result is the $n=2$ case of \cite[Theorem 2.5]{BEJ_ndimWA}.
Note also that in the introduction, and in the statement of Theorem~\ref{t:headline}, we restricted ourselves to cocycles on $A$ taking values in~$A^*$. This is no loss of generality: by the $n=2$ case of \cite[Corollary 2.11]{BEJ_ndimWA}, if there exists some $X$ with $\cZ^2_{\rm alt}(A,X) \neq0$ then $\simpaltco{2}{A}\neq 0$.

The following example illustrates what the abstract definitions above mean in practice, and provides motivation for later constructions.

\begin{eg}[Derivations and cocycles on $C^1(\bbT)$ and $C^1(\bbT^2)$]
\label{eg:C1-cases}
\
\begin{romnum}

\item
Define a linear map $D: C^1(\bbT)\to C^1(\bbT)^*$ by
\[ D(f_1)(f_0) = \int_{\bbT} \frac{\partial f_1}{\partial\theta}(p)f_0 (p)\,dp \]
where $p=e^{2\pi i\theta}$ and $dp$ denotes the usual uniform measure on $\bbT$. It follows immediately from the product rule that $D$ is a derivation.
\item
Define a bilinear map $F: C^1(\bbT^2)\times C^1(\bbT^2)\to C^1(\bbT^2)^*$ by
\[
F(f_1,f_2)(f_0) = \int_{\bbT^2} \left(
\frac{\partial f_1}{\partial\theta_1}(p)
\frac{\partial f_2}{\partial\theta_2}(p)
-
\frac{\partial f_1}{\partial\theta_2}(p)
\frac{\partial f_2}{\partial\theta_1}(p)
\right) f_0(p) \,dp \]
where $p=(e^{2\pi i\theta_1},e^{2\pi i\theta_2})$ and $dp$ denotes the usual uniform measure on $\bbT^2$.
Clearly $F$ is alternating as a bilinear map, and a similar calculation to part (i) shows it is a derivation in the first variable; thus it is an alternating $2$-cocycle.
\end{romnum}
\end{eg}

\begin{rem}
We noted earlier that an alternating $2$-cocycle is a coboundary if and only if it is identically zero. This is true for all $n\geq 1$ by \cite[Prop.~2.9]{BEJ_ndimWA}, and so there is a natural \emph{injection} of the vector space $\cZ^n_{\rm alt}(A,X)$ into the $n$th Hochschild cohomology group $\cH^n(A,X)$. The range of this injection is one summand in a canonical decomposition of $\cH^n(A,X)$ into $n$~pieces; for further discussion of this decomposition, see e.g.~\cite[\S3]{YC_QJM10} and \cite[\S9.4]{Weibel}.
\end{rem}

\begin{rem}\label{r:HKR behind the scenes}
We briefly leave the world of Banach algebras to mention some important context from algebraic geometry.
If ${\mathsf R}$ is the coordinate ring of a smooth complex variety, every cocycle on ${\sf R}$ with symmetric coefficients is equivalent to an alternating one; this is one version of the \dt{Hochschild--Kostant--Rosenberg theorem.} For instance, this applies to the complex coordinate ring of the algebraic group $\SL_2$, which occurs naturally as a dense subalgebra of $\FA(\SU(2))$. While the HKR theorem itself does not seem to extend to the Banach-algebraic setting, it suggests that the alternating cocycles on commutative Banach algebras have some deeper meaning, rather than being ad hoc definitions, and hence deserve further study.
\end{rem}

\end{subsection}

\begin{subsection}{Tools for constructing alternating cocycles}
\label{ss:Z-alt-tools}
Following a strategy analogous to those used in \cite{CG_WAAG1,CG_WAAG2,LLSS_WAAG}, we shall prove Theorem~\ref{t:headline} by establishing the non-vanishing of $\simpaltco{2}{\FA(G_1)}$ for some judiciously chosen closed subgroup $G_1\subset G$.
We record the following lemma for later reference.

\begin{lem}\label{l:pull back}
Let $\cA$ and $\cB$ be commutative Banach algebras and let $\theta:\cA\to \cB$ be a continuous homomorphism. Then for any $F\in\simpaltco{2}{\cB}$, the induced map $\theta^*F$ defined by
\[
\theta^* F (a_1,a_2)(a_0) \defeq F(\theta(a_1),\theta(a_2))(\theta(a_0))
\]
belongs to $\simpaltco{2}{\cA}$. If $F\neq 0$ and $\theta$ has dense range, then $\theta^*F\neq 0$.
\end{lem}

The proof follows easily from the definitions and we omit the details.

As mentioned in Section~\ref{ss:intro-justify}, in commutative algebra there is a standard procedure for constructing alternating $2$-cocycles on a tensor product of two algebras, given a pair of derivations on the respective algebras.
With minor modifications, one can do the same in the setting of commutative Banach algebras and symmetric Banach bimodules. This observation is surely known to specialists in the cohomology of Banach algebras, but we have not found an explicit statement in the literature; this is somewhat surprising since it provides a natural converse to a special case of \cite[Theorem 3.6]{BEJ_ndimWA}.

\begin{lem}[possibly folklore]\label{l:Z2alt-on-bsptp}
Given commutative Banach algebras $A$ and $B$, symmetric Banach bimodules $X$ and $Y$ over $A$ and $B$ respectively, and derivations $D_A:A\to X$, $D_B:B\to Y$, the formula
\begin{equation}\label{eq:BA-wedge}
F( a_1\tp b_1, a_2\tp b_2) \defeq
 [D_A(a_1)\cdot a_2 ] \tp [b_1\cdot D_B(b_2)]
 - [a_1\cdot D_A(a_2) \tp D_B(b_1)\cdot b_2]
\end{equation}
 defines an alternating $2$-cocycle
 $F : A\ptp_\gamma B \times A\ptp_\gamma B \to X\ptp_\gamma Y$.
 \end{lem}

Since some of the relevant calculations will recur when we come to prove Theorem~\ref{t:mainthm}, we provide a detailed proof.

 \begin{proof}
 Since $D_A$ and $D_B$ are bounded, $F$ extends to a continuous bilinear map
 \newline
 $(A\bsptp B)\times (A\bsptp B)\to X\bsptp Y$, by the universal property of $\bsptp$. Given $a_0,a_1,a_2\in A$ and $b_0,b_1,b_2\in B$, direct calculation yields
\[
 F(a_1\tp b_1, a_2\tp b_2) + F(a_2\tp b_2, a_1\tp b_1) = 0
\]
(using the fact that $X$ and $Y$ are \emph{symmetric} bimodules)
 and
\[ \begin{aligned}
 F(a_1a_2\tp b_1b_2, a_0\tp b_0)
& =
\left\{ \begin{gathered}
{[D_A(a_1a_2)\cdot a_0 ] \tp [b_1b_2\cdot D_B(b_0)]} \\
- [a_1a_2\cdot D_A(a_0) \tp D_B(b_1b_2)\cdot b_0]
\end{gathered} \right.
\\
& =
\left\{ \begin{gathered}
{[D_A(a_2)\cdot a_1a_0 ] \tp [b_1b_2\cdot D_B(b_0)]} \\ 
+ {[D_A(a_1)\cdot a_2a_0 ] \tp [b_1b_2\cdot D_B(b_0)]} \\
- [a_1a_2\cdot D_A(a_0) \tp D_B(b_2)\cdot b_1b_0] \\
- [a_1a_2\cdot D_A(a_0) \tp D_B(b_1)\cdot b_2b_0] \\
\end{gathered} \right.
\\
& = F(a_2\tp b_2 , a_1a_0\tp b_1b_0) + F(a_1\tp b_1, a_2a_0\tp b_2b_0)
\end{aligned} \]
(using the fact that $D_A$ and $D_B$ are derivations into symmetric bimodules).
Hence, by linearity and continuity, $F$ is both alternating and a derivation in the first variable, so by the earlier remarks in this section it is an alternating $2$-cocycle as required.
\end{proof}
 
The formula \eqref{eq:BA-wedge} should be compared with Example~\ref{eg:C1-cases}.
In that example, $F$ was obtained by ``wedging together'' two derivations defined on the dense subalgebra $C^1(\bbT)\tp C^1(\bbT)$, one being a copy of $D$ acting in the $\theta_1$ direction and the other being a copy of $D$ acting in the $\theta_2$ direction. Lemma~\ref{l:Z2alt-on-bsptp} may be regarded as an abstract analogue of this construction.
However, since $C^1(\bbT)\bsptp C^1(\bbT) \neq C^1(\bbT^2)$, the lemma does not suffice on its own to construct alternating $2$-cocycles on $C^1(\bbT^2)$.
Indeed, the next example shows that for some function algebras on $\bbT^2$, the approach suggested by Lemma~\ref{l:Z2alt-on-bsptp} cannot possibly work.

\begin{eg}\label{eg:lip-alg}\
For $\alpha\in (0,1)$ and $n\in \Nat$,
consider the \dt{little Lipschitz algebra} $\lip_\alpha(\bbT^n)\equiv \lip(\bbT^n,d(\cdot)^\alpha)$. Classical Fourier analysis tells us that the trigonometric polynomials are dense in $\lip_\alpha(\bbT^n)$.
Moreover, given $f_1$ and $f_2$ in $\lip_\alpha(\bbT)$, a straightforward calculation shows that $f_1\tp f_2\in \lip_\alpha(\bbT^2)$. Thus we may identify the algebraic\footnotemark\ tensor product $\lip_\alpha(\bbT)\tp \lip_\alpha(\bbT)$ with a dense subalgebra $\sR\subset \lip_\alpha(\bbT^2)$.
\footnotetext{This embedding extends to a continous map $\phi:\lip_\alpha(\bbT)\bsptp\lip_\alpha(\bbT)\to \lip_\alpha(\bbT^2)$. In fact, $\phi$ is injective, but the proof requires some technical facts from Banach space theory which would distract us here. }

\begin{figure}[htp]
\hfil\begingroup
\renewcommand*{\arraystretch}{1.6}
\begin{tabular}{l|c|c|c}
       & $0< \alpha \leq 1/2$ & $1/2 <\alpha \leq 2/3$ & $2/3\leq \alpha < 1$ \\
  \hline
  $\simpco{1}{\lip_\alpha(\bbT)}$ & zero & non-zero & non-zero\\
  $\simpaltco{2}{\lip_\alpha(\bbT)}$ & zero & zero & zero \\
  \hline
  $\simpco{1}{\lip_\alpha(\bbT^2)}$ & zero & non-zero & non-zero \\
  $\simpaltco{2}{\lip_\alpha(\bbT^2)}$ & zero & zero & non-zero \\
\end{tabular}
\endgroup
\hfil

\caption{$1$-cocycles and $2$-cocycles for $\lip_\alpha(\bbT)$ and $\lip_\alpha(\bbT^2)$}
\label{fig:lip-info}
\end{figure}

We now appeal to some consequences of \cite[Corollary 4.4]{BEJ_ndimWA}; the relevant information is displayed in Figure~\ref{fig:lip-info}.
When $1/2< \alpha <1$, $\lip_\alpha(\bbT)$ has a non-zero derivation into its dual. Using Lemma~\ref{l:Z2alt-on-bsptp}, we obtain an $F\in \cZ^2(\lip_\alpha(\bbT) \bsptp \lip_\alpha(\bbT), \lip_\alpha(\bbT)^* \bsptp \lip_\alpha(\bbT)^*)$, which by inspection is non-zero when restricted to~$\sR$.

Since $\lip_\alpha(\bbT)^* \bsptp \lip_\alpha(\bbT)^*$ embeds continuously in $\lip_\alpha(\bbT^2)^*$, we therefore have a natural, non-zero, densely-defined alternating $2$-cocycle on $\lip_\alpha(\bbT^2)$ for all $\alpha\in (1/2, 1)$. However, when $\alpha \in (1/2, 2/3]$, this cannot be extended to give an  element of $\simpaltco{2}{\lip_\alpha(\bbT^2)}$, since the latter space vanishes.
\end{eg}

Lemma~\ref{l:Z2alt-on-bsptp} can still be applied to produce alternating $2$-cocycles on $\FA(H)\bsptp\FA(L)$. However, this is not enough to produce cocycles on $\FA(H\times L)$, because of the following two facts.
\begin{enumerate}
\renewcommand{\labelenumi}{\arabic{enumi})}
\item The natural map $\FA(H)\bsptp \FA(L) \to \FA(H\times L)$ is surjective if and only if either $H$ or $L$ has an abelian subgroup of finite index \cite{Los_tp-AG}. (See also \cite[\S3.6]{KL_AGBGbook}, with the caveat that they write $\ptp$ instead of $\bsptp$.)
\item If $H$ has an abelian subgroup of finite index, then the only (continuous) derivation $\FA(H)\to\FA(H)^*$ is the zero map. (In fact it suffices that the connected component of $H$ be abelian; see \cite[Theorem 3.3]{ForRun_amenAG} or \cite[\S4.5]{KL_AGBGbook}.)
\end{enumerate}

\end{subsection}
\end{section}

\begin{section}{Co-cb maps and Fourier algebras}
\label{s:UNNAMED}

This section is devoted to the infrastructure needed for the proof of Theorem~\ref{t:mainthm}. We pay particular attention to issues of functoriality; the reason for introducing operator space tensor products and the opposite operator space structure is not just to equip Banach spaces with extra structure, but to be able to combine linear maps that respect this extra structure.

\begin{subsection}{Operator spaces, tensor products, and co-cb maps}\label{ss:OS-prelim}
All concepts not defined explicitly here can be found in standard sources, such as the early chapters of \cite{ER_OSbook} or \cite{Pis_OSbook}.

Henceforth, we abbreviate the phrase ``operator space structure'' to \oss.
Given operator spaces $X$ and $Y$, $\CB(X,Y)$ denotes the space of completely bounded maps $X\to Y$; note that this space has a canonical \oss, defined via the identification $M_n\CB(X,Y)\iso \CB(X,M_nY)$.

Whenever $H$ is a Hilbert space and we refer to $B(H)$ as an operator space, we assume (unless explicitly stated otherwise) that it is equipped with its usual, canonical \oss; note that if we do this, then there is a natural and completely isometric identification of $B(H)$ with $\CB(\COL_H)$, where $\COL_H$ denotes $H$ equipped with the column \oss.

\para{Tensor products and tensor norms}
The projective and injective tensor products of operator spaces 
are denoted by $\ptp$ and $\itp$ respectively (this is the notation of \cite{ER_OSbook}, rather than that of \cite{Pis_OSbook}).
Note that if $E$ and $F$ are Hilbert spaces then the underlying Banach space of $\COL_E \ptp (\COL_F)^*$ is $E\bsptp F^*$, but the underlying Banach space of $\COL_E \ptp \COL_F$ is $E\tp_2 F$.

Given operator spaces $V$, $W$ and $X$, we say that a bilinear map $V\times W \to X$ is \dt{jointly completely bounded}\footnotemark\ (j.c.b.) if it extends to a completely bounded map $V\ptp W \to X$.
\footnotetext{This seems to now be the accepted terminology, and agrees with \cite{Pis_OSbook}. One should beware that in \cite{ER_OSbook} such bilinear maps are called ``completely bounded'', which in most later sources is instead used for those maps linearized by the Haagerup tensor product.}
This is equivalent to saying that the ``curried map'' $V\to L(W,X)$ extends to a completely bounded map $V\to \CB(W,X)$, or the same with $V$ and $W$ interchanged.
Indeed $V\ptp W$ may be characterized, up to completely isometric isomorphism, as the completion of $V\tp W$ that satisfies
\begin{equation}\label{eq:cb-curry}
\CB(V\ptp W,X) \iso \CB(V,\CB(W,X)) \iso \CB(W,\CB(V,X)) \text{ completely isometrically.}
\tag{$\diamondsuit$}
\end{equation}

If $f\in \CB(E,X)$ and $g\in \CB(F,Y)$ then by tensoring we obtain completely bounded maps $E\ptp F\to X\ptp Y$ and $E\itp F\to X\itp Y$; for extra emphasis, these maps will be denoted by $f\ptp g$ and $f\itp g$ respectively.

We shall make passing use of the Haagerup tensor norm, but only ``at level~1'', and we only require the following facts:
\begin{enumerate}
\renewcommand{\labelenumi}{\arabic{enumi})}
\renewcommand{\theenumi}{\arabic{enumi})}
\item\label{li:htp-fact1}
 for any $w\in E\tp F$ we have $\norm{w}_{E\smallitp F} \leq \norm{w}_{E\htp F} \leq \norm{w}_{E\ptp F}$;
\item\label{li:htp-fact2}
if $A$ and $B$ are $\Cst$-algebras and $w\in A\tp B$, then
\[
\norm{w}_{A\htp B} = \inf_{n\in\Nat} \inf \left\{
\left\Vert\sum\nolimits_{j=1}^n a_j^*a_j\right\Vert^{1/2}
\left\Vert\sum\nolimits_{j=1}^n b_jb_j^*\right\Vert^{1/2}
\right\}
\]
where the inner infimum is over all representations of $w$ as $\sum_{j=1}^n a_j\tp b_j$.
\end{enumerate}
For a proof of \ref{li:htp-fact1} see e.g. \cite[Theorem 9.2.1]{ER_OSbook};
for \ref{li:htp-fact2}, see e.g. \cite[Ch.~5]{Pis_OSbook}.

\para{The opposite operator space structure}
Given an operator space $W$, one may define a new sequence of matrix norms on $W$ by
\[ \Norm{ \sum\nolimits_i a_i \tp w_i}_{(n), {\rm opp}} \defeq \norm{ \sum\nolimits_i a_i^\top \tp w_i}_{(n)} \qquad(a_i\in \Mat_n, w_i\in W).\]
(C.f.~\cite[\S2.10]{Pis_OSbook}.)
These satisfy Ruan's axioms and hence equipe $W$ with a new \oss, which we call the \dt{opposite \oss}; the resulting operator space will be denoted by $\til{W}$.
 In longer expressions, when considering the opposite operator space, we use the notation {$\longopp{(\dots)}$}; for instance $\longopp{B(H)}$ denotes $B(H)$ equipped with the opposite of its usual \oss.

Note that for a Hilbert space $H$, $\longopp{(\COL_H)}=\ROW_H$ and $\longopp{(\ROW_H)}= \COL_H$.
More generally, $\longopp{(\til{W})}=W$.

\begin{rem}
$\til{W}$ is often denoted in the literature by $W^{\rm op}$. We have chosen different notation because there is a potential conflict with the usage of $A^{\rm op}$ to denote the ``opposite algebra'', i.e. the algebra with the same underlying vector space but with reversed product. The two conventions match happily if $A=B(H)$ but are at odds if $A=\FA(G)$.
\end{rem}

It is easily checked that if $f:X \to Y$ is completely bounded, then so is $f: \til{X} \to \til{Y}$, with the same cb-norm. To emphasise the functorial behaviour we write this as $\til{f}: \til{X}\to\til{Y}$.
The same calculation gives, with some book-keeping, a more precise result: we omit the details.

\begin{lem}\label{l:functorial-opp}
Given operator spaces $X$ and $Y$, the assignment $f\mapsto \til{f}$ defines a completely isometric isomorphism $\longopp{\CB(X,Y)} \cong \CB(\til{X},\til{Y})$. In particular, we can identify $\longopp{(X^*)}$ with~$(\til{X})^*$.
\end{lem}

Note that for any operator spaces $V$ and $W$, the identity map on $V\tp W$ extends to a completely isometric isomorphism $\til{V}\ptp\til{W} \cong \longopp{(V\ptp W)}$. One can show this using the explicit definition of the matrix norms associated to~$\ptp$, but it can also be deduced from the characterization in~\eqref{eq:cb-curry}, combined with repeated application of Lemma~\ref{l:functorial-opp}.

\paragraph{Co-cb maps between operator spaces.}
The next definition is non-standard (though it has some precedent\footnotemark\ in \cite{Pis_co-cb}) but will be extremely useful for statements and calculations later on.

\begin{dfn}[Co-complete boundedness]
Let $V$ and $W$ be operator spaces and let $f:V\to W$ be a linear map.
Note that $f$ is c.b.\ from $V$ to $\til{W}$ if and only if it is c.b. from $\til{V}$ to~$W$; in either case we say that $f:V\to W$ is \dt{co-completely bounded} (\dt{co-cb} for short).
Similarly,  $f$ is a complete isometry from $V$ to $\til{W}$ if and only if it is a complete isometry from $\til{V}$ to $W$; we then say that $f:V\to W$ is a \dt{co-complete isometry.}
\end{dfn}
\footnotetext{In \cite{Pis_co-cb} this concept is called ``completely co-bounded'', but then abbreviated to ``co-cb'' just as we have done. Our terminology is chosen to avoid potential confusion with cohomology theory (``cobounded'') or operator theory (``coisometry'').}

The notion of co-completely bounded map seems to have gone largely unmentioned or unstudied in the literature. One notable exception is \cite{Pis_co-cb}, which sets up some general machinery and obtains interesting results in connection with Schur multipliers.

\begin{lem}\label{l:adjoint flips oss}
Let $F$ be a Hilbert space. 
The $\Cplx$-linear map $B(F)\to B(F^*)$ that sends an operator $b\in B(F)$ to its Banach-space adjoint $b^\# : F^*\to F^*$ is a co-complete isometry.
\end{lem}

{\bf Warning:} as we will see in the proof, our chosen notational conventions are important: $B(F^*)$ is given the \oss\ of $\CB(\COL_{F^*})$ rather than $\CB((\COL_F)^{*})=\CB(\ROW_{F^*})$.

\begin{proof}
For any operator spaces $X$ and $Y$, the map $\CB(X,Y)\to \CB(Y^*,X^*)$ defined by taking adjoints is a \emph{complete}\footnotemark\ isometry.
\footnotetext{This is standard knowledge but we could not locate an explicit statement of this result in the literature. It follows easily from the completely isometric identifications $\CB(Y^*,X^*)\cong\CB(Y^*\ptp X, \Cplx)\cong\CB(X,Y^{**})$, since the map $\CB(X,Y)\to \CB(Y^*,X^*)$ then corresponds to the canonical inclusion $\CB(X,Y)\to\CB(X,Y^{**})$.}
Thus $b\mapsto b^\#$ defines a complete isometry
\[
B(F)\equiv \CB(\COL_F) \to \CB((\COL_F)^*) = \CB(\ROW_{F^*}) = \CB(\longopp{(\COL_{F^*})}).
\]
Taking opposites and applying Lemma~\ref{l:functorial-opp}, the result follows.
\end{proof}

\begin{rem}[Transpose in a basis-free setting]
\label{r:basis-free-transpose}
Consider $L^2(\Om)$ for some measure space $(\Om,\mu)$ (we suppress mention of the measure $\mu$ for notational convenience), and let $\alpha: \eta\mapsto \evetabar$ denote the canonical, \emph{linear}, isometric isomorphism $L^2(\Om)\to L^2(\Om)^*$
that is defined via Equation~\eqref{eq:L^2-iso-to-dual}.
This defines a normal $*$-isomorphism ${\rm Ad}_{\alpha}:B(L^2(\Om)^*) \to B(L^2(\Om))$.
Calculation shows that composing ${\rm Ad}_\alpha$ with the map $\#: B(L^2(\Om)) \to B(L^2(\Om)^*)$ from Lemma~\ref{l:adjoint flips oss} yields a \emph{linear}, co-complete isometry
\[ \top: B(L^2(\Om) \to B(L^2(\Om)) \quad,\quad b\mapsto b^\top \;, \]
where $b^\top(\xi) \defeq \overline{ b^*\overline{\xi}}$ for each $\xi\in L^2(\Om)$.

Note that if $\Om$ is countable and infinite, equipped with counting measure, then $\top$ is just the usual ``transpose'' operator for an infinite matrix.
\end{rem}
\end{subsection}

\begin{subsection}{The operator space structure on a Fourier algebra}\label{ss:AG-oss}
The results in this section are all known to specialists, but are included here for the reader's convenience, and to ensure that we have consistent notation and conventions.

Given a Hilbert space $H$ and a ``concrete'' von Neumann algebra $\cM\subseteq B(H)$: the predual $\cM_*$ is a natural quotient of $B(H)_*= \COL_H \ptp (\COL_H)^*$, and may thus be equipped with  the quotient \oss; moreover, if we take the dual of this \oss\ on $\cM_*$, we recover the original subspace \oss\ on $\cM$.
(See \cite[Prop.\ 4.2.2]{ER_OSbook}.)
Thus $(\cM_*)^*\cong\cM$ completely isometrically.

In particular, consider the \dt{group von Neumann algebra} $\VN(G)\subseteq B(L^2(G))$. Then one can identify $\FA(G)$ with the unique isometric predual of $\VN(G)$. Our chosen definition for $\FA(G)$ allows one to deduce this quickly by considering the adjoint of the map in Equation~\eqref{eq:left-reg coeff},
\[ \Psi^*: \FA(G)^* \to (L^2(G)\bsptp L^2(G)^*)^* \cong B(L^2(G)), \]
observing that for each $s\in G$, $\Psi^*$ maps the character ${\rm ev}_s$ to the translation operator $\lambda(s)$.
See \cite[Lemma~2.8.2]{KL_AGBGbook} for details.
We now take, as our canonical \oss\ on $\FA(G)$, the one induced by $\Psi$ from $\VN(G)_*$.

\begin{rem}[Tomato, tomato]
\label{r:tomato-tomato}
Care is needed when combining parts of the literature. Some sources, following the general framework of locally compact quantum groups, define $\FA(G)$ to be the subspace of $C_0(G)$ obtained by identifying a vector functional $\omega_{\xi,\eta} \in \VN(G)_*$ with the function $s\mapsto \ip{\lambda(s^{-1})\xi}{\eta}_{L^2(G)}$.
While this gives the same Banach function algebra $\FA(G)\subset C_0(G)$ as in this paper, it yields the \emph{opposite} \oss\ to the one we have just defined.
This is related to Proposition~\ref{p:known-AG-os}\ref{li:check-is-cci} below.
\end{rem}

Given $f\in\Cplx^G$, let $\AGcheck{f}(x)=f(x^{-1})$. For any $\xi,\eta\in L^2(G)$ and $x\in G$ we have
\[
\Psi(\xi\tp\evetabar)(x^{-1}) =
\int_G \xi(xs)\eta(s)\,ds =
\int_G \eta(x^{-1}s)\xi(s)\,ds =
\Psi(\eta\tp \evxibar)(x)
\,.
\]
It follows that the map $f\mapsto \AGcheck{f}$ defines a contractive involution on $\FA(G)$ (which must therefore be isometric).
This is known as the \dt{flip map} or \dt{check map} on $\FA(G)$.

\begin{rem}\label{r:transpose of check}
Direct calculations show that when we identify $\FA(G)^*$ with $\VN(G)$, the adjoint of the check map on $\FA(G)$ coincides with the restriction to $\VN(G)$ of the transpose operator $\top: B(L^2(G))\to B(L^2(G))$.
This observation is folklore (and a similar calculation may be found in the proof of \cite[Prop.~1.5(ii)]{ForRun_amenAG}).
\end{rem}

\begin{prop}[Known results needed later]
\label{p:known-AG-os}
Let $G$, $G_1$ and $G_2$ be locally compact groups.
\begin{alphnum}
\item\label{li:AG-ptp}
The natural inclusion $\FA(G_1)\tp \FA(G_2) \to \FA(G_1\times G_2)$ extends to a completely isometric isomorphism $\FA(G_1)\ptp\FA(G_2) \cong \FA(G_1\times G_2)$.
\item\label{li:VN-itp}
The natural inclusion $\VN(G_1)\tp \VN(G_2)\to \VN(G_1\times G_2)$ extends to a complete isometry $\VN(G_1)\itp\VN(G_2) \to \VN(G_1\times G_2)$, whose range is the minimal $\Cst$-tensor product $\VN(G_1)\mintp \VN(G_2)$.
\item\label{li:check-cb-iff-VA}
The check map on $\FA(G)$ is completely bounded if and only if $G$ has an (open) abelian subgroup of finite index.
\item\label{li:check-is-cci}
The check map on $\FA(G)$ is a co-complete isometry.
\end{alphnum}
\end{prop}

\begin{proof}
Parts \ref{li:AG-ptp} and \ref{li:VN-itp} are general results about von Neumann algebra preduals and $\Cst$-algebras; see e.g.~\cite[Theorem 7.2.4]{ER_OSbook} and \cite[Prop.~8.1.6]{ER_OSbook}) respectively.
Part \ref{li:check-cb-iff-VA} is stated as \cite[Prop.~1.5(ii)]{ForRun_amenAG}; the authors give a complete proof, while making it clear that the result was already known to previous specialists.
Part \ref{li:check-is-cci} holds because the adjoint of the check map coincides with the transpose operator (see Remark~\ref{r:transpose of check}), which is a co-complete isometry by Remark~\ref{r:basis-free-transpose}.
\end{proof}

\end{subsection}
\end{section}

\begin{section}{Constructing $2$-cocycles from co-cb derivations}
\label{s:make cocycles}

We start in some generality, since the preliminary results may be useful in subsequent work.

The following terminology is not entirely standard, but is analogous to the more familiar notions of \textit{completely contractive Banach algebra} and \textit{completely contractive Banach (bi)module} that have appeared in the literature.
By a \dt{cb-Banach algebra}, we mean an operator space $A$ equipped with a bilinear, j.c.b.~and associative map $A\times A \to A$. Given such an $A$, we define a \dt{cb-Banach $A$-bimodule} to be an operator space $X$, equipped with an $A$-bimodule structure such that the left action $A\times X \to X$ and the right action $X\times A \to X$ are both~j.c.b.

Clearly $A$ itself is a cb-Banach $A$-bimodule; it is also routine to check that if $X$ is a cb-Banach $A$-bimodule, so is $X^*$ when equipped with the dual \oss.
These notions also interact well with the ``opposite \oss\ functor''. If $A$ is a cb-Banach algebra then so is $\til{A}$; and if $X$ is a cb-Banach $A$-bimodule, $\til{X}$ is a cb-Banach $\til{A}$-bimodule.

\begin{rem}\label{r:uh oh}
Given a cb-Banach algebra $A$, the class of cb-Banach $A$-bimodules is usually not closed under the operation of taking opposites. For instance, suppose $A$ is unital. If $\til{A}$ were a cb-Banach $A$-bimodule, then for each $x\in A$ the orbit map $a\mapsto ax$ would be completely bounded from $A$ to $\til{A}$. Taking $x=1_A$ we conclude that the identity map on $A$ is co-cb. In particular, if $A=\FA(G)$ for $G$ compact, this would force $G$ to be virtually abelian (combine parts \ref{li:check-cb-iff-VA} and \ref{li:check-is-cci} of Proposition~\ref{p:known-AG-os}).
\end{rem}

\begin{prop}\label{p:twist-wedge}
Let $A$ and $B$ be cb-Banach algebras; let $X$ be a cb-Banach $A$-bimodule and $Y$ a cb-Banach $B$-bimodule. Let $T_A\in \CB(A, \til{X})$ and $T_B\in \CB(B, \til{Y})$. Then, if we define $F_1, F_2: (A\tp B)\times (A\tp B) \to X\tp Y$ by
\[ \begin{aligned}
F_1(a_1\tp b_1, a_2\tp b_2) & = [T_A(a_1)\cdot a_2] \tp [b_1\cdot T_B(b_2)]\;, \\
F_2(a_1\tp b_1, a_2\tp b_2) & = [a_1\cdot T_A(a_2)]\tp [T_B(b_1)\cdot b_2 ]\;,\\
\end{aligned}
\]
both $F_1$ and $F_2$ extend to bounded bilinear maps $(A\ptp B)\times (A\ptp B) \to X\ptp \til{Y}$.
\end{prop}

\begin{proof}
We will only give the proof for $F_1$; the proof for $F_2$ is very similar.

Since $T_A:A \to \til{X}$ is completely bounded, so is $\til{T_A}: \til{A} \to X$. Therefore, if we put
\[
S(a_1\tp b_1\tp a_2\tp b_2) \defeq T_A(a_1) \tp  b_1 \tp a_2 \tp T_B(b_2)
\]
we obtain a complete contraction
\[
S = \til{T_A}\ptp \iota_{\til{B}} \ptp \iota_A \ptp T_B :
 \til{A}\ptp \til{B} \ptp A\ptp B \longrightarrow
 X \ptp \til{B} \ptp A \ptp \til{Y} \;.
\]
Also, since $X$ is a cb-Banach $A$-bimodule  and $\til{Y}$ is a cb-Banach $\til{B}$-bimodule,
putting
\[
R(x\tp b_1\tp a_2\tp y) = (x\cdot a_2)\tp (b_1\cdot y)
\]
defines a complete contraction $R: X \ptp \til{B} \ptp A \ptp \til{Y}\to X\ptp \til{Y}$.

Since $\til{A}\ptp \til{B} = (A\ptp B)^\sim$, the composite map $RS$ defines a j.c.b.~bilinear map from $(A\ptp B)^\sim \times (A\ptp B)$ to $X\ptp\til{Y}$, which agrees with $F_1$ on $(A\tp B) \times (A\tp B)$. In particular, $F_1$ extends to a bounded bilinear map (no longer completely bounded!) from $(A\ptp B) \times (A\ptp B)$ to $X\ptp\til{Y}$.
\end{proof}

\begin{rem}
The proof of Proposition~\ref{p:twist-wedge} would have been much easier if $T_A\tp \iota_B \tp \iota_A \tp T_B$ extended to a continuous linear map from $(A\ptp B)\ptp_\gamma (A\ptp B)$ to $(X\ptp A)\ptp_\gamma (B\ptp Y)$.
However, we see no reason why this should always hold,
as it requires an interchange/distributivity result for $\ptp$ and~$\bsptp$.
\end{rem}

In view of the earlier formula \eqref{eq:BA-wedge}, one would like to apply Proposition~\ref{p:twist-wedge} with $T_A$ and $T_B$ being co-cb derivations into symmetric cb-bimodules. However, this stops short of producing genuine $2$-cocycles: the resulting bilinear map merely takes values in $X\ptp \til{Y}$,
and in view of Remark~\ref{r:uh oh} there is no reason to suppose that this is even a Banach $A\ptp B$-bimodule. (It is a Banach $A\bsptp B$-bimodule, but that does not help us.)
To go further, we need to move from $X\ptp \til{Y}$ to $X\itp Y$, and this is where we require Theorem~\ref{t:twisted inclusion}, whose proof we now turn to.
For convenience we recall the statement of the theorem.

\paragraph{\thf Theorem \ref{t:twisted inclusion} ({\itshape reprise}).}
{\itshape 
Let $X$ and $Y$ be operator spaces. Then
$\norm{w}_{X\smallitp\til{Y}} \leq \norm{w}_{X\ptp Y}$ for all $w\in X\tp Y$,
and so the identity map on $X\tp Y$ extends to a contraction $\theta_{X,Y}: X\ptp Y\to X\itp \til{Y}$.
}

\begin{proof}[Proof of Theorem~\ref{t:twisted inclusion}]
We start by reducing to a special case.

\para{Step 1}
Given  a pair of operator spaces $X$ and $Y$, 
let $j_X:X\to B(E)$ and $j_Y:Y\to B(F)$ be completely isometric embeddings, for some Hilbert spaces $E$ and $F$.
Note that $\til{j_Y} :\til{Y} \to \longopp{B(Y)}$ is also a complete isometry.

Suppose we know Theorem~\ref{t:twisted inclusion} holds for the particular operator spaces $B(E)$ and $B(F)$.
Then we have a diagram as shown in Figure \ref{fig:embedding trick},
in which the left-hand vertical arrow is a (complete) contraction, while the right-hand vertical arrow is a (complete) \emph{isometry} (since $\itp$ respects complete isometries).

\begin{figure}[htp]
\[ \begin{diagram}[tight,width=7em,height=3em]
B(E) \ptp \longopp{B(F)} & \rTo^{\theta_{B(E),B(F)}}  & B(E)\itp B(F)  \\
\uTo^{j_X\ptp \til{j_Y}} & & \uTo_{j_X\itp j_Y} \\
X \ptp \til{Y}  & \rDots & X\itp Y
\end{diagram} \]
\caption{An embedding trick}
\label{fig:embedding trick}
\end{figure}

Moreover, the diagram in Figure~\ref{fig:embedding trick} ``commutes on elementary tensors''.
Hence, for any $z\in X\tp Y$, we have
\[
\norm{z}_{X\smallitp \til{Y}}
=
\norm{(j_X \tp j_Y)(z)}_{B(E)\smallitp \longopp{B(F)}}
\leq
\norm{(j_X\tp j_Y)(z)}_{B(E)\ptp B(F)}
\leq
\norm{z}_{X\ptp Y}
\;.
\]

Thus: {\itshape if Theorem~\ref{t:twisted inclusion} holds for $B(E)$ and $B(F)$, then it holds for all operator spaces.}

\para{Step 2}
Observe that if $E$ and $F$ are Hilbert spaces and $w\in B(E)\tp B(F)$, the norm of $w$ in $B(E)\itp \longopp{B(F)}$ coincides with the norm of the associated elementary operator on $S_2(F,E)$, the space of Hilbert-Schmidt operators $F\to E$. This is a variation on a well-known fact in $\Cst$-algebra theory that can be found in various sources; to avoid any notational ambiguity,  we give a precise statement in the following lemma.

\begin{lem}
Define $\Phi_2: B(E)\tp B(F) \to B(S_2(F,E))$ by
$\Phi_2(a\tp b)(c) = acb$.
Then for all $w\in B(E)\tp B(F)$ we have
\begin{equation}\label{eq:re-express}
\norm{\Phi_2(w)}_{B(S_2(F,E))}
= \norm{w}_{B(E)\smallitp \longopp{B(F)}}
\,.
\end{equation}
\end{lem}

\begin{proof}
Recall that if $b\in B(F)$ then $b^\#:B(F^*)\to B(F^*)$ denotes its Banach-space adjoint. Now we make two observations.
Firstly: by Lemma~\ref{l:adjoint flips oss} and the fact $\itp$ respects complete isometries, $\iota\tp\#$ is a complete isometry from $B(E)\itp \longopp{B(F)}$ onto $B(E)\itp B(F^*)$.
Secondly: there is an injective $*$-homomorphism $\theta:B(E)\mintp B(F^*)\to B(S_2(F,E))$, which on elementary tensors satisfies
$\theta(a\tp b^\#)(c) =acb$.
Since $\Phi_2=\theta\circ(\iota\tp\#)$, Equation~\eqref{eq:re-express} holds.

(To justify the second point in a little more detail: let $\alpha:E\basetp_2 F^* \to S_2(F,E)$ be the Hilbert-space isomorphism which sends $x\tp\phi$ to $y\mapsto \phi(y)x$; then $\alpha$ intertwines the natural $*$-representation of the incomplete algebra $B(E)\tp B(F^*)$ on $E\basetp_2 F^*$ with the map $\theta$. See also \cite[Prop.~2.9.1]{Pis_OSbook} or the calculations preceding \cite[Eqn.~(3.5.1)]{ER_OSbook}.)
\end{proof}

\para{Step 3}
Combining Steps 1 and 2, we see that Theorem~\ref{t:twisted inclusion} will follow if we can prove the following claim:
\textit{given $E$, $F$ and $\Phi_2$ as in Step 2, the function $\Phi_2$ extends to a contractive linear map $B(E)\ptp B(F) \to B(S_2(F,E))$.}
(We remind the reader that if $E$ and $F$ are infinite-dimensional, one cannot expect this map to be completely bounded.)

Our proof of the claim is based on an interpolation argument.
For $p\in [1,\infty]$ let $\SCH_p(F,E)$ denote the space of Schatten-$p$ operators from $F\to E$, equipped with its standard norm.
We adopt the convention that $\SCH_\infty(F,E) =K(F,E)$, the space of all compact operators $F\to E$, equipped with the operator norm. 
Define $\Phi_p: B(E)\tp B(F) \to B(\SCH_p(F,E))$ by the formula
\[
\Phi_p(a\tp b)(c) = acb
\qquad\qquad(\text{$a\in B(E)$, $c\in \SCH_p(F,E)$, $b\in B(F)$}).
\]
When $p=2$ this is consistent with our earlier notation.

\begin{lem}\label{l:rainwater}
Let $E$ and $F$ be Hilbert spaces, and let $w\in B(E)\tp B(F)$. Let $\sigma : B(E)\tp B(F) \to B(F)\tp B(E)$ denote the flip map $x\tp y\mapsto y\tp x$.
\begin{romnum}

\item\label{li:Sinf}
$\norm{\Phi_{\infty}(w) : \SCH_\infty(F,E)\to \SCH_\infty(F,E)} \leq \norm{w}_{B(E)\htp B(F)}$.
\item\label{li:S1}
$\norm{\Phi_1(w) : \SCH_1(F,E)\to \SCH_1(F,E)} \leq \norm{\sigma(w)}_{B(F)\htp B(E)}$.
\item\label{li:S2}
$\norm{\Phi_2(w) : \SCH_2(F,E)\to \SCH_2(F,E)} \leq \left( \norm{w}_{B(E)\htp B(F)} \norm{\sigma(w)}_{B(F)\htp B(E)} \right)^{1/2}$.
\end{romnum}
\end{lem}

\begin{proof}
Part \ref{li:Sinf} follows immediately by quoting Haagerup's theorem that $\Phi_\infty$ extends to a (complete) isometry from $B(E)\htp B(F)$ to $\CB(\SCH_\infty(F,E))$.
However, there is also a direct easy proof: given $w=\sum_{j=1}^n a_j\tp b_j$, it suffices to show that
\[
\Norm{\sum\nolimits_{j=1}^n a_j cb_j} \leq \norm{c} \Norm{\sum\nolimits_{j=1}^n a_ja_j^*}^{1/2} \Norm{\sum\nolimits_{j=1}^n b_j^*b_j}^{1/2} \;.
\]
This follows from standard calculations with ``row'' and ``column'' block matrices: for details, see e.g.\cite[Remark 1.13]{Pis_OSbook}, in particular the formula (1.12) in \cite{Pis_OSbook}.

Part \ref{li:S1} follows from part \ref{li:Sinf} and duality. In more detail: given $w= \sum_{j=1}^n a_j \tp b_j \in B(E)\tp B(F)$, consider the elementary operator $\Phi_\infty(\sigma(w))$ defined on $\SCH_\infty(E,F)$ by  $d\mapsto \sum_{j=1}^nb_jda_j$.
By part \ref{li:Sinf}, applied with the roles of $E$ and $F$ reversed, $\norm{\Phi_\infty(\sigma(w))} \leq \norm{\sigma(w)}_h$. On the other hand, consider the standard trace pairing between $\SCH_1(F,E)$ and $\SCH_\infty(E,F)$, where $s \in \SCH_1(F,E)$ acts as the functional $t\mapsto \Tr(st)$.
Straightforward calculations show that with respect to this pairing,
the Banach-space adjoint of $\Phi_\infty(\sigma(w)): \SCH_\infty(E,F)\to \SCH_\infty(E,F)$ is the elementary operator $\Phi_1(w): \SCH_1(F,E)\to \SCH_1(F,E)$. Since a linear map and its adjoint have the same norm, \ref{li:S1} is proved.

Finally, note that the elementary operator defined by $w$ acts simultaneously on all $\SCH_p(F,E)$ for $p\in [1,\infty]$. Since $\SCH_\infty$ and $\SCH_1$ form an interpolation couple with $(\SCH_1,\SCH_\infty)_{1/2}=\SCH_2$, part~\ref{li:S2} now follows from parts~\ref{li:Sinf} and~\ref{li:S1} by applying the Riesz--Thorin interpolation theorem.
\end{proof}

To finish off, note that the os-projective tensor norm dominates both the Haagerup tensor norm and the ``reversed'' Haagerup tensor norm. More precisely:
for arbitrary operator spaces $V_1$ and $V_2$ and $w\in V_1\tp V_2$, we have
\[
\norm{w}_{V_1\htp V_2} \leq \norm{w}_{V_1\ptp V_2}
\quad\text{and}\quad
\norm{\sigma(w)}_{V_2\htp V_1} \leq \norm{\sigma(w)}_{V_2\ptp V_1}=
\norm{w}_{V_1\ptp V_2} \;.
\]
Combining these inequalities with Lemma~\ref{l:rainwater}\ref{li:S2}, we have verified the claim at the beginning of Step 3, and this completes the proof of Theorem~\ref{t:twisted inclusion}.
\end{proof}

\begin{rem}\label{r:could say more}
We can strengthen the conclusion of Theorem \ref{t:twisted inclusion}, although at present we do not have applications of the stronger version.
Going back to Step 1, we can run the same argument with the os-projective tensor product replaced by either the Haagerup tensor product or its reversed version. Combining this with Step 2 and Lemma~\ref{l:rainwater}\ref{li:S2}, we conclude that for arbitrary operator spaces $X$ and $Y$ and $w\in X\tp Y$,
\begin{equation}
\norm{w}_{X\smallitp \til{Y}}
\leq 
\left( \norm{w}_{X\htp Y} \norm{\sigma(w)}_{Y\htp X} \right)^{1/2}
\leq 
\frac{1}{2}\left( \norm{w}_{X\htp Y} + \norm{\sigma(w)}_{Y\htp X} \right)
.
\end{equation}
\end{rem}

We now have the necessary ingredients for our main technical theorem.

\begin{thm}\label{t:mainthm}
Let $H$ and $L$ be locally compact groups. Suppose that there exist non-zero, co-cb derivations $D_H: \FA(H)\to \FA(H)^*$ and $D_L:\FA(L)\to\FA(L)^*$.
Then the bilinear map $F_0: (\FA(H)\tp \FA(L) )\times (\FA(H)\tp \FA(L) ) \longrightarrow \FA(H)^* \tp \FA(L)^*$ defined by
\[
F_0( a_1\tp b_1, a_2\tp b_2) \defeq
 [D_A(a_1)\cdot a_2 ] \tp [b_1\cdot D_B(b_2)]
 - [a_1\cdot D_A(a_2) \tp D_B(b_1)\cdot b_2]
 \]
extends to a non-zero alternating $2$-cocycle $F:\FA(H\times L) \times \FA(H\times L) \to \FA(H\times L)^*$.

Consequently, for any locally compact group $G$ which contains a closed isomorphic copy of $H\times L$, we have $\simpaltco{2}{\FA(G)}\neq 0$.
\end{thm}

\begin{proof}
For this proof only, just to ease notation slightly, we denote the Fourier algebras of $G$, $H$, $L$ and $H\times L$ by $A_G$, $A_H$, $A_L$ and $A_{H\times L}$ respectively; it is convenient to sometimes denote their duals  by $V_G$, $V_H$, $V_L$ and $V_{H\times L}$ respectively.

Let $D_H: A_H \to V_H$ and $D_L: A_L\to V_L^*$ be non-zero, co-cb derivations, and let
\[ F_0: (A_H\tp A_L )\times (A_H\tp A_L ) \to V_H \tp V_L \]
be as defined in the statement of the theorem.
By Proposition~\ref{p:twist-wedge}, $F_0$ extends to a bounded bilinear map from
$(A_H\ptp A_L )\times (A_H\ptp A_L )$ to $V_H \ptp \longopp{V_L}$. Applying Theorem~\ref{t:twisted inclusion} with $X= V_H$ and $Y=\longopp{V_L}$, we obtain a bounded bilinear map
\[ F: (A_H\ptp A_L )\times (A_H\ptp A_L ) \to V_H \itp V_L \]
that extends $F_0$.
Recall (Proposition~\ref{p:known-AG-os}) that the natural map $A_H\ptp A_L\to A_{H\times L}$ is a (completely isometric) isomorphism and that $V_H \itp V_L$ embeds (completely isometrically) in $V_{H\times L}=(A_{H\times L})^*$. Thus $F$ can be viewed as a bilinear map $A_{H\times L} \times A_{H\times L } \to (A_{H\times L})^*$.

As in the proof of Lemma~\ref{l:Z2alt-on-bsptp}, the defining formula for $F_0$ shows that $F$ is an alternating $2$-derivation on the dense subalgebra $A_H\tp A_L\subset A_{H\times L}$. By the usual continuity argument we deduce that $F\in \simpaltco{2}{A_{H\times L}}$.

We now show that $F$ is not identically zero. Since $F_0$ takes values in $V_H \tp V_L$ and the natural map $V_H\tp V_L \to V_{H\times L}$ is injective, it suffices to show that $F_0$ is not identically zero. Observe that if $a\in \FA_H$ and $b\in \FA_L$ we have
\begin{equation*}\label{eq:hacky}
\begin{aligned}
F_0(a^3\tp b, a\tp b) & = [D_H(a^3)\cdot a] \tp [b\cdot D_L(b)] - [a^3\cdot D_H(a)] \tp [D_L(b)\cdot b] \\
& = 2a^3\cdot D_H(a) \tp b\cdot D_L(b) \\
& = \frac{1}{4} D_H(a^4)\tp D_L(b^2).
\end{aligned}
\end{equation*}
By Lemma~\ref{l:dense-powers}, elements of the form $a^4$ span a dense subspace of $A_H$, and elements of the form $b^2$ span a dense subspace of $A_L$.
Therefore, since $D_H$ is continuous and non-zero, there exists $a\in A_H$ such that $D_H(a^4)\neq 0$; similarly, there exists $b\in A_L$ such that $D_L(b^2)\neq 0$. We conclude that $F_0(a^3\tp b, a\tp b)\neq 0$, as required.

This proves the first part of the theorem. The second part follows by pulling back the non-zero $2$-cocycle $F\in \simpaltco{2}{A_{H\times L}}$ along the restriction homomorphism $A_G \to A_{H\times L}$ (see Lemma \ref{l:pull back}).
\end{proof}

To use Theorem~\ref{t:mainthm} effectively, we need to know examples of $H$ for which such a $D_H$ exists. It turns out that the very first non-zero derivation constructed from a Fourier algebra to its dual, which was produced by Johnson in \cite{BEJ_AG}, can be shown with hindsight to be~co-cb!
In fact, during the writing of \cite{CG_WAAG1},  the present author and Ghandehari had already observed that if $H$ is one of the groups
\begin{alphnum}
\item\label{li:CG_SU(2)}
$\SU(2)$ or $\SO(3)$;
\item\label{li:CG_ax+b}
the real $ax+b$ group (the connected component of $\Real\rtimes\Real^*$);
\item\label{li:CG_Hred}
the reduced Heisenberg group (the quotient of the $3$-dimensional real Heisenberg group by a central copy of $\Zahl$);
\end{alphnum}
then in each case, the explicit non-zero derivation $D_H: \FA(H)\to\FA(H)^*$ that is described in \cite{CG_WAAG1} turns out to be co-cb. Showing this requires some work, but is mostly just a matter of composing $D_H$ with (the adjoint of) the check map and using the Plancherel theorem for each group, c.f.~the formulas and remarks in \cite[\S7]{CG_WAAG1}.
However, this observation was never written down, since at the time we did not have any applications of~it.

Let us see how these facts, combined with Theorem~\ref{t:mainthm}, yield the first two cases of Theorem~\ref{t:headline}. For $n\geq 4$, we may embed $\GL_2(\Cplx)\times \GL_2(\Cplx)$ as a closed subgroup of $\GL_n(\Cplx)$ by sending $(g_1,g_2)$ to the block-diagonal matrix ${\rm diag}(g_1,g_2, I_{n-4})$; the same construction works with $\Cplx$ replaced by $\Real$. If $H_n$ denotes any of $\SU(n)$, $\SL(n,\Real)$ or ${\rm Isom}(\Real^n)$, then our embedding maps $H_2\times H_2$ onto a closed subgroup of $H_n$. Note also that the real $ax+b$ group is isomorphic to the standard parabolic subgroup of $\SL(2,\Real)$. Therefore, we may combine Theorem~\ref{t:mainthm} with the examples \ref{li:CG_SU(2)} and \ref{li:ax+b} mentioned above.

To obtain the remaining case of Theorem~\ref{t:headline}, it would suffice to exhibit a non-zero co-cb derivation from $\FA({\rm Isom}(\Real^2))$ to its dual.
For this, the results of \cite{CG_WAAG1,CG_WAAG2} are insufficient and we require results from the subsequent paper \cite{LLSS_WAAG}.
In fact, one can use results from that paper to obtain alternative proofs for the cases \ref{li:CG_SU(2)} \ref{li:CG_ax+b} and \ref{li:CG_Hred}, and so for clarity of exposition we devote the next section to summarizing and making use of the relevant parts of \cite{LLSS_WAAG}.

\begin{rem}
An alternative proof that the Fourier algebra of the real $ax+b$ group supports a non-zero co-cb derivation, independent of both \cite{CG_WAAG1} and \cite{LLSS_WAAG}, will appear as part of the forthcoming work~\cite{CG_affdualconv}.
\end{rem}

\end{section}

\begin{section}{Obtaining co-cb derivations}
\label{s:co-cb-der-AG}

In this section, we show in Theorem~\ref{t:underyournoses} that there is a plentiful supply of non-zero co-cb derivations from Fourier algebras to their duals.
For the strongest results in this direction, we make use of the hard work done by the authors of \cite{LLSS_WAAG} in proving the Lie case of the Forrest--Runde conjecture.
While Theorem~\ref{t:underyournoses}  is not hard to invent if one reads \cite{LLSS_WAAG} in its entirety, it is never actually stated in that paper. We shall therefore extract some of the components which are used to prove \cite[Theorem~3.2]{LLSS_WAAG}, and reassemble them into a ``black box'' that will be more suitable for our purposes.

For $G$ a Lie group, let $C_c^\infty(G)$ denote the space of compactly supported smooth functions on~$G$. This is contained in $\FA(G)$ by \cite[(3.26)]{Eym_BSMF64} and is easily seen to be dense in $\FA(G)$ (since $C_c^\infty(G)$ is dense in $L^2(G)$ and is closed under convolution).

\begin{prop}[Lee--Ludwig--Samei--Spronk, \cite{LLSS_WAAG}]\label{p:portmanteau}
Let $H$ be any one of the following (connected, real) Lie groups:
\begin{alphnum}
\item\label{li:SU(2)}
$\SU(2)$ or $\SO(3)$;
\item\label{li:ax+b}
the real $ax+b$ group (the connected component of $\Real\rtimes\Real^*$);
\item\label{li:Hred}
the reduced Heisenberg group (the quotient of the $3$-dimensional real Heisenberg group by a central copy of $\Zahl$);
\item\label{li:E(2)}
the Euclidean motion group of $\Real^2$;
\item\label{li:grelaud}
the ``Gr\'elaud groups'' $G_\theta$ (certain semidirect products $\Real^2\rtimes_\theta \Real$ where $\theta$ parametrizes the eigenvalues of the corresponding action of $\Real$ on the Lie algebra of $\Real^2$)
\end{alphnum}
Then there exist a weight function $v\in L^1(H)$, not identically zero, and an element $X$ of the Lie algebra of $H$, such that when we take the corresponding Lie derivative $\lieder{X}:C_c^\infty(H)\to C_c^\infty(H)$,
\begin{equation}
\label{eq:twisted-der}
\left\vert \int_H (\lieder{X}\tp\iota)u(s,s^{-1})v(s)\,ds \right\vert \leq \norm{u}_{\FA(H\times H)}
\qquad\text{for all $u\in C_c^\infty(H\times H)$.}
\end{equation}
\end{prop}

\begin{proof}
In each case, there is a calculation in \cite{LLSS_WAAG} that provides suitable $X$ and $v$. (Strictly speaking, $X$ and $v$ are chosen together with $S=S_{X,v}\in\VN(H\times H)$ such that the integral in \eqref{eq:twisted-der} agrees with $\langle S, u\rangle$ for all $u\in C_c^\infty(H\times H)$. By density and continuity arguments, if such an $S$ exists it is uniquely determined, and by rescaling the weight function $v$ we can always arrange that $\norm{S}\leq1$.)

For \ref{li:SU(2)}, see \cite[Theorem 2.4]{LLSS_WAAG} --- strictly speaking, the cited result only proves this for $\SU(2)$, but the same calculation using representation theory and orthogonality relations goes through for $\SO(3)$, almost word for word.

For \ref{li:ax+b}, see \cite[Theorem~2.9]{LLSS_WAAG};
for \ref{li:Hred}, see \cite[Theorem~2.6]{LLSS_WAAG};
for \ref{li:E(2)}, see \cite[Theorem~2.5]{LLSS_WAAG};
and
for \ref{li:grelaud}, see \cite[Theorem~2.8]{LLSS_WAAG}.
\end{proof}

\begin{rem}\label{r:precis of LLSS}
In \cite{LLSS_WAAG} the results assembled in Proposition~\ref{p:portmanteau} were used as follows.
Recall that a commutative Banach algebra $A$ is said to be \dt{weakly amenable} if $\cZ^1(A,A^*)=0$. 
Let $H$ be one of the groups in Proposition~\ref{p:portmanteau}: then by \cite[Lemma~2.1]{LLSS_WAAG}, the functional on $\FA(H\times H)$ that is uniquely defined by \eqref{eq:twisted-der} serves as a witness that $H$ has the following property:
\begin{equation}
\tag{AD}
\text{the anti-diagonal in $H\times H$ is not a set of local smooth synthesis.}
\end{equation}
By \cite[Theorem 1.6]{LLSS_WAAG}, the universal cover of a connected Lie group $L$ has property (AD) if and only if $L$ does; and by \cite[Theorem 1.3]{LLSS_WAAG}, if $L$ has property (AD) then $\FA(L)$ is not weakly amenable. Now by the structure theory of real Lie algebras, every non-abelian connected Lie group $G$ contains a closed Lie subgroup $H_0$ with the same Lie algebra as one of the groups listed in Proposition~\ref{p:portmanteau}; it follows that $\FA(H_0)$, and hence $\FA(G)$, fails to be weakly amenable, as predicted by the Forrest--Runde conjecture.
\end{rem}

As indicated by the previous remark, the authors of \cite{LLSS_WAAG} did not pursue Proposition~\ref{p:portmanteau} with the goal of constructing explicit derivations on Fourier algebras, as they aimed to establish stronger structural properties for a wider class of groups.
For the present paper, what matters is the following consequence of Proposition~\ref{p:portmanteau}, which appears to be a new observation.

\begin{thm}\label{t:underyournoses}
Let $H$ be one of the groups listed in {\rm(a)--(e)} of Proposition~\ref{p:portmanteau}, and let $X$ and $v$ be as provided by that proposition.
Then there is a unique bounded linear map $D:\FA(H)\to\FA(H)^*$ that satisfies
\begin{equation}\label{eq:underyournoses}
D(f_1)(f_0) = \int_H (\lieder{X} f_1)(s) f_0(s) v(s)\,ds \qquad\text{for all $f_0,f_1\in C_c^\infty(H)$.}
\end{equation}
Furthermore, $D$ is a non-zero co-cb derivation.
\end{thm}

\begin{proof}
First, observe that there is at most one (norm-)continuous function $D:\FA(H)\to\FA(H)^*$ satisfying \eqref{eq:underyournoses}, because $C_c^\infty(H)$ is dense in~$\FA(H)$.

Since $X$ and $v$ are chosen so that the inequality \eqref{eq:twisted-der} is satisfied, there exists a unique $\psi\in \FA(H\times H)^*$ satisfying
\begin{equation}\label{eq:en route}
\begin{aligned}
 \psi(f \tp g)
& = \int_H (\lieder{X} \tp \iota)(f\tp g)(s,s^{-1})v(s)\,ds  \\
& = \int_H (\lieder{X} f)(s)\ g(s^{-1})\  v(s)\,ds
&\quad\text{for all $f,g\in C_c^\infty(H)$.}
\end{aligned} \tag{$*$}
\end{equation}
Let $T:\FA(H)\to\FA(H)^*$ be the completely bounded map corresponding to $\psi$, and define $D:\FA(H)\to\FA(H)^*$ by
$D(f_1)(f_0) \defeq T(f_1)(\AGcheck{f_0})$ for $f_0,f_1\in\FA(H)$.
Since the check map is a co-complete isometry (Proposition~\ref{p:known-AG-os}\ref{li:check-is-cci}), $D:\FA(H)\to\FA(H)^*$ is co-cb. Moreover, since $\psi$ satisfies \eqref{eq:en route}, $D$ satisfies \eqref{eq:underyournoses}.

It remains to show that $D$ is a derivation. Since $C_c^\infty(H)$ is dense in $\FA(H)$ and $D$ is norm-continuous, it suffices to check that
$D(g_1g_2)(g_0)=D(g_2)(g_0g_1)+D(g_1)(g_2g_0)$ for all $g_0,g_1,g_2\in C_c^\infty(H)$. This follows from \eqref{eq:underyournoses} and the fact that $\lieder{X}:C_c^\infty(H)\to C_c^\infty(H)$ is a derivation.
\end{proof}

Theorem~\ref{t:headline} now follows by combining Theorem~\ref{t:mainthm} with the cases \ref{li:SU(2)}, \ref{li:ax+b} and \ref{li:E(2)} of Theorem \ref{t:underyournoses},
following the argument described at the end of Section~\ref{s:make cocycles}.

\begin{rem}[Infinite-dimensional spaces of derivations]
Inspecting the results used to prove Proposition~\ref{p:portmanteau}, one sees that for each of the solvable cases, and an appropriately chosen $X$, the set of $v$ which ``work'' contains an infinite-dimensional vector space. It follows that in Theorem \ref{t:underyournoses} one obtains not just a single $D$ of the desired form, but an infinite-dimensional\footnotemark\ vector space of such derivations.%
\footnotetext{This also fits with what one expects from considering derivations on more general Banach function algebras~$A$; often there exists a dense subalgebra $\sA_0\subset A$ such that the natural predual of $\cH^1(A,A^*)$ contains a free $\sA_0$-module in the algebraic sense, forcing $\cH^1(A,A^*)$ to be infinite-dimensional as a vector space.}
Thus for these groups, the suggestion at the end of \cite[\S2.4]{LLSS_WAAG},
concerning the space of bounded derivations from the Fourier algebra to its dual, cannot be correct as stated.
\end{rem}
\end{section}

\begin{section}{Avenues for further work}\label{s:finalremarks}
In future work, we intend to set out a more systematic study of the higher-degree alternating cocycles on Fourier algebras, with the intention of exploring an associated numerical invariant that can be viewed as a kind of ``dimension'' associated to such algebras.
Since one would like to calcuate or estimate this numerical invariant for as many small examples as possible, progress on the following natural question could be a useful guide for future work.

\begin{qn}
Does Theorem~\ref{t:headline} remain true for $n=2$ or $n=3$?
\end{qn}

Currently our guess is that the answer is negative for $\SU(2)$ and $\SL(2,\Real)$, and positive for $\SU(3)$, $\SL(3,\Real)$, ${\rm Isom}(\Real^2)$ and ${\rm Isom}(\Real^3)$, but there is insufficient evidence to support any firm conjectures at this stage.

Turning to Theorem~\ref{t:twisted inclusion}: one would like to understand better the comparison map $\theta_{X,Y} : X\ptp Y \to X\itp \til{Y}$, perhaps by making greater use of the sharper result outlined in Remark~\ref{r:could say more}. Indeed, a natural next step is to repeat the (complex) interpolation argument used in Lemma~\ref{l:rainwater} at the level of operator spaces and cb-norms of elementary operators, to see what $\theta_{X,Y}$ looks like at higher matrix levels.

The co-cb derivations that are crucial to proving Theorem~\ref{t:headline} provide natural examples of c.b.\ maps from $\longopp{\FA(G)}$ to $\VN(G)$ that behave like noncommutative Fourier multipliers (this is not immediately apparent from what is stated in Section~\ref{s:co-cb-der-AG}, but can be seen by inspecting the details in \cite{CG_WAAG1} and \cite{LLSS_WAAG}).

\begin{qn}
Given that $\longopp{\FA(G)}$ and $\VN(G)$ are the endpoints of the scale of noncommutative $L^p$-spaces associated to $\VN(G)$, are there other Fourier multipliers from $L^p(\VN(G)) \to L^r(\VN(G))$ which satisfy some form of the Leibniz identity?
For fixed $p$ and $r$, what can we say about the space of such multipliers?
\end{qn}

We finish with some natural questions concerning co-cb derivations on Fourier algebras, which are all aimed at strengthening or sharpening the conclusion of Theorem~\ref{t:underyournoses}.

\begin{qn}
The derivations constructed in \cite{CG_WAAG1} for $\SU(2)$, the real $ax+b$ group and the reduced Heisenberg group are all cyclic and co-cb (c.f.~the construction in \cite{CG_affdualconv}). Is every cyclic derivation on a Fourier algebra automatically co-cb?
\end{qn}

\begin{qn}\label{q:co-cb_CG2}
Let $G$ be the 3-dimensional real Heisenberg group. The results of \cite{CG_WAAG2} construct a non-zero derivation $D$ from $\FA(G)$ to a certain symmetric Banach $\FA(G)$-bimodule $\sW$. Can $\sW$ be made into a cb-Banach $\FA(G)$-bimodule in such a way that $D:\FA(G)\to\sW$ is co-cb?
\end{qn}

\begin{qn}\label{q:co-cb_LLSS}
It was shown in \cite{LLSS_WAAG} that the property (AD) for a Lie group $G$, mentioned in Remark~\ref{r:precis of LLSS}, ensures that there is a non-zero derivation $D:\FA(G)\to\FA(G)^*$. Does it also guarantee that one can choose $D$ to be co-cb?
\end{qn}

In view of the good hereditary properties of (AD), a positive answer to Q\ref{q:co-cb_LLSS} would allow us to transfer co-completely bounded derivations between Fourier algebras of Lie groups which have the same universal cover, and hence by using the strategy outlined in Remark~\ref{r:precis of LLSS} one could strengthen Theorem~\ref{t:underyournoses} to the following result: {\itshape every non-abelian connected Lie group $H$ has non-zero co-cb derivations from $\FA(H)$ to $\FA(H)^*$.}

\begin{qn}\label{q:co-cb_Losert}
Can the explicit derivations constructed by Losert \cite{Los_WAAG_pre} on connected groups that are not necessarily Lie, be made into co-cb derivations from Fourier algebras into cb-Banach bimodules?
\end{qn}
The constructions in \cite{Los_WAAG_pre} are closer in spirit to \cite{CG_WAAG2} than to \cite{LLSS_WAAG}, and so Question~\ref{q:co-cb_CG2} would serve as a warm-up for Question~\ref{q:co-cb_Losert}.

\end{section}

\section*{Acknowledgments}

A preliminary announcement of some of these results, with a different emphasis and less comprehensive results, was circulated in 2016 as an unpublished preprint; the author thanks M. Daws and A. Skalski for several comments and corrections on that document.
He also thanks the authors of \cite{LLSS_WAAG} for useful discussions about technical aspects of their paper, and to V. Losert for sharing a copy of the preprint~\cite{Los_WAAG_pre}.

Less direct, but no less important, thanks are due to
M. Ghandehari and E. Samei for their interest and encouragement over several years concerning alternating cocycles on Fourier algebras, and to M. Whittaker for a conversation at the British Mathematical Colloquium 2019 in Lancaster which prompted the author to finally write up this work as a proper paper.

Some of these results were presented in conference talks at the Abstract Harmonic Analysis meeting in Kaohsiung, 2018, and the International Workshop on Harmonic Analyis and Operator Theory, Istanbul, 2019. The author thanks the organizers of these meetings for their respective invitations to present this work, and for bringing together speakers from a broad range of specialist interests for enjoyable discussions.

Most of the writing of this article was done during the COVID-19 pandemic, under a period of lockdown conditions in England. The author would therefore like to thank various colleagues at Lancaster University for an online mixture of  camaraderie, commiseration, complaints, and computer support during these months, which has gone some way towards replicating a normal working environment. He hopes to one day visit Barnard Castle.


\vfill

\noindent
Department of Mathematics and Statistics\\
Fylde College, Lancaster University\\
Lancaster, United Kingdom LA1 4YF

\smallskip
\noindent
Email: \texttt{y.choi1@lancaster.ac.uk}

\end{document}

%% file: zaltagostp_mac.tex
\usepackage{amsmath,amssymb,amsfonts}

\usepackage{sectsty}
\allsectionsfont{\sffamily}

\numberwithin{equation}{section}

\renewcommand{\emph}[1]{{\sl#1}\/}
\newcommand{\defeq}{\mathbin{:=}}
\newcommand{\dt}[1]{{\itshape#1}\/}

\newcommand{\iso}{\cong}

\newcommand{\im}{\operatorname{im}}

\newcommand{\para}[1]{\paragraph{#1.}} 

\newcommand{\norm}[1]{{\Vert{#1}\Vert}}
\newcommand{\Norm}[1]{{\bigl\Vert{#1}\bigr\Vert}}

\newcommand{\basetp}{{\otimes}}
\newcommand{\tp}{\mathbin{\basetp}}
\newcommand{\ptp}{\mathbin{\widehat{\basetp}}}
\newcommand{\bsptp}{\mathbin{\widehat{\basetp}}_\gamma}
\newcommand{\mintp}{\mathbin{\basetp}_{\rm min}}
\newcommand{\htp}{\mathbin{\basetp_h}}

\newcommand{\itp}{\mathbin{\basetp^\vee}}
\newcommand{\smallitp}{\basetp^\vee}

\newcommand{\AGcheck}[1]{{#1}^{\raise0.4ex\hbox{$\scriptscriptstyle\vee$}}}

\newcommand{\ip}[2]{{\langle#1,#2\rangle}}
\newcommand{\evetabar}{{\rm ev}_{\overline{\eta}}}
\newcommand{\evxibar}{{\rm ev}_{\overline{\xi}}}

\newcommand{\Tr}{\operatorname{\rm Tr}}

\makeatletter
\newcommand{\optfs}{\@ifnextchar.{}{.\@}}
\makeatother
\newcommand{\oss}{o.s.s\optfs}

\newcommand{\COL}{\textsf{COL}}
\newcommand{\ROW}{\textsf{ROW}}

\newcommand{\SCH}{{\mathcal S}} 
\renewcommand{\SCH}{{S}} 

\newcommand{\Om}{\Omega}

\newcommand{\cA}{{\mathcal A}}
\newcommand{\cB}{{\mathcal B}}
\newcommand{\cC}{{\mathcal C}}
\newcommand{\cF}{{\mathcal F}}
\newcommand{\cH}{{\mathcal H}}
\newcommand{\cM}{{\mathcal M}}

\newcommand{\cU}{{\mathcal U}}
\newcommand{\cZ}{{\mathcal Z}}

\newcommand{\mbb}[1]{{\mathbb#1}}

\newcommand{\Cplx}{\mbb{C}}
\newcommand{\Real}{\mbb{R}}
\newcommand{\Zahl}{\mbb{Z}}
\newcommand{\Nat}{\mbb{N}}
\newcommand{\Mat}{\mbb{M}}

\newcommand{\bbT}{\mbb{T}}

\newcommand{\CB}{\operatorname{\it CB}} 

\newcommand{\til}[1]{\widetilde{#1}}
\newcommand{\longopp}[1]{{#1}^{\sim}} 

\newcommand{\lin}{\operatorname{lin}}

\newcommand{\lip}{\operatorname{lip}\nolimits}

\newcommand{\simpco}[2]{{\mathcal Z}^{#1}(#2, {#2}^*)}
\newcommand{\simpaltco}[2]{{\mathcal Z}_{\rm alt}^{#1}(#2, {#2}^*)}

\newcommand{\SU}{{\rm SU}}
\newcommand{\SO}{{\rm SO}}
\newcommand{\SL}{{\rm SL}}
\newcommand{\GL}{{\rm GL}}

\newcommand{\lieder}[1]{{\partial}_{{#1}}}  

\newcommand{\sA}{{\sf A}}
\newcommand{\sB}{{\sf B}}
\newcommand{\sR}{{\sf R}}
\newcommand{\sW}{{\sf W}}
\newcommand{\sX}{{\sf X}}
\newcommand{\sY}{{\sf Y}}

\newcommand{\VN}{\operatorname{VN}}
\newcommand{\FA}{\operatorname{A}}

\newcommand{\Cst}{\operatorname{C}^*}

\newenvironment{romnum}{%
\begin{enumerate}
\setlength{\itemsep}{0em}
\renewcommand{\labelenumi}{{\rm(\roman{enumi})}}
\renewcommand{\theenumi}{{\rm(\roman{enumi})}}
}{\end{enumerate}\ignorespacesafterend}

\newenvironment{alphnum}{%
\begin{enumerate}
\setlength{\itemsep}{0em}
\renewcommand{\labelenumi}{{\rm(\alph{enumi})}}
\renewcommand{\theenumi}{{\rm(\alph{enumi})}}
}{\end{enumerate}\ignorespacesafterend}

\RequirePackage[PS,nohug]{diagrams}

\newarrow{Sub}       C--->  
\newarrow{Quot}    ----{>>} 
\newarrow{Eq} ===== 
\newarrow{Dash} {}{dash}{}{dash}>
\newarrow{Dots} ....>

\RequirePackage{amsthm}   

\newcounter{pulse}[section]
\numberwithin{pulse}{section}  

\newcommand{\thf}{\sc} 

\theoremstyle{plain}
\newtheorem{thm}[pulse]{\thf Theorem}

\newtheorem{prop}[pulse]{\thf Proposition}
\newtheorem{lem}[pulse]{\thf Lemma}

\theoremstyle{definition}
\newtheorem{dfn}[pulse]{\thf Definition}

\newtheorem{eg}[pulse]{\thf Example}
\theoremstyle{remark}
\newtheorem{rem}[pulse]{\thf Remark}

\newcounter{qnum}
\newtheorem{qn}[qnum]{\thf Question}

%% file: widehack.tex
\DeclareFontFamily{U}{mathx}{\hyphenchar\font45}
\DeclareFontShape{U}{mathx}{m}{n}{
      <5> <6> <7> <8> <9> <10>
      <10.95> <12> <14.4> <17.28> <20.74> <24.88>
      mathx10
      }{}
\DeclareSymbolFont{mathx}{U}{mathx}{m}{n}
\DeclareFontSubstitution{U}{mathx}{m}{n}
\DeclareMathAccent{\widecheck}{0}{mathx}{"71}

\renewcommand{\itp}{\widecheck{\basetp}}
\renewcommand{\smallitp}{\widecheck{\basetp}}
\renewcommand{\itp}{\widecheck{\basetp}}


%% file: Z2alt-AG-ostp_arX4.bbl
\begin{thebibliography}{10}

\bibitem{YC_QJM10}
{\sc Y.~Choi}, {\em Hochschild homology and cohomology of {$\ell^1({\bf
  Z}^k_+)$}}, Q. J. Math., 61 (2010), pp.~1--28.

\bibitem{CG_WAAG1}
{\sc Y.~Choi and M.~Ghandehari}, {\em Weak and cyclic amenability for {F}ourier
  algebras of connected {L}ie groups}, J. Funct. Anal., 266 (2014),
  pp.~6501--6530.

\bibitem{CG_WAAG2}
\leavevmode\vrule height 2pt depth -1.6pt width 23pt, {\em Weak amenability for
  {F}ourier algebras of 1-connected nilpotent {L}ie groups}, J. Funct. Anal.,
  268 (2015), pp.~2440--2463.

\bibitem{CG_affdualconv}
\leavevmode\vrule height 2pt depth -1.6pt width 23pt, {\em Dual convolution for
  the affine group of the real line}, Complex Anal. Oper. Theory, to appear (2021).
\newblock 
See arXiv \href{https://arxiv.org/abs/2009.05497}{2009.05497}


\bibitem{ER_OSbook}
{\sc E.~G. Effros and Z.-J. Ruan}, {\em Operator spaces}, vol.~23 of London
  Mathematical Society Monographs. New Series, Oxford University Press, New
  York, 2000.

\bibitem{Eym_BSMF64}
{\sc P.~Eymard}, {\em L'alg\`ebre de {F}ourier d'un groupe localement compact},
  Bull. Soc. Math. France, 92 (1964), pp.~181--236.

\bibitem{ForRun_amenAG}
{\sc B.~E. Forrest and V.~Runde}, {\em Amenability and weak amenability of the
  {F}ourier algebra}, Math. Z., 250 (2005), pp.~731--744.

\bibitem{HK_AP-group}
{\sc U.~Haagerup and J.~Kraus}, {\em Approximation properties for group
  {$C^*$}-algebras and group von {N}eumann algebras}, Trans. Amer. Math. Soc.,
  344 (1994), pp.~667--699.

\bibitem{BEJ_AG}
{\sc B.~E. Johnson}, {\em Non-amenability of the {F}ourier algebra of a compact
  group}, J. London Math. Soc. (2), 50 (1994), pp.~361--374.

\bibitem{BEJ_ndimWA}
\leavevmode\vrule height 2pt depth -1.6pt width 23pt, {\em Higher-dimensional
  weak amenability}, Studia Math., 123 (1997), pp.~117--134.

\bibitem{KL_AGBGbook}
{\sc E.~Kaniuth and A.~T.-M. Lau}, {\em Fourier and {F}ourier-{S}tieltjes
  algebras on locally compact groups}, vol.~231 of Mathematical Surveys and
  Monographs, American Mathematical Society, Providence, RI, 2018.

\bibitem{LLSS_WAAG}
{\sc H.~H. Lee, J.~Ludwig, E.~Samei, and N.~Spronk}, {\em Weak amenability of
  {F}ourier algebras and local synthesis of the anti-diagonal}, Adv. Math., 292
  (2016), pp.~11--41.

\bibitem{Los_tp-AG}
{\sc V.~Losert}, {\em On tensor products of {F}ourier algebras}, Arch. Math.
  (Basel), 43 (1984), pp.~370--372.

\bibitem{Los_WAAG_pre}
\leavevmode\vrule height 2pt depth -1.6pt width 23pt, {\em On weak amenability
  of {F}ourier algebras}.
\newblock Preprint, 2019.

\bibitem{Pis_OSbook}
{\sc G.~Pisier}, {\em Introduction to operator space theory}, vol.~294 of
  London Mathematical Society Lecture Note Series, Cambridge University Press,
  Cambridge, 2003.

\bibitem{Pis_co-cb}
\leavevmode\vrule height 2pt depth -1.6pt width 23pt, {\em Completely
  co-bounded {S}chur multipliers}, Oper. Matrices, 6 (2012), pp.~263--270.

\bibitem{Sam_cb-loc-der}
{\sc E.~Samei}, {\em Bounded and completely bounded local derivations from
  certain commutative semisimple {B}anach algebras}, Proc. Amer. Math. Soc.,
  133 (2005), pp.~229--238.

\bibitem{Spr_OWA}
{\sc N.~Spronk}, {\em Operator weak amenability of the {F}ourier algebra},
  Proc. Amer. Math. Soc., 130 (2002), pp.~3609--3617.

\bibitem{Takesaki_vol2}
{\sc M.~Takesaki}, {\em Theory of operator algebras. {II}}, vol.~125 of
  Encyclopaedia of Mathematical Sciences, Springer-Verlag, Berlin, 2003.
\newblock Operator Algebras and Non-commutative Geometry, 6.

\bibitem{Weibel}
{\sc C.~A. Weibel}, {\em An introduction to homological algebra}, vol.~38 of
  Cambridge Studies in Advanced Mathematics, Cambridge University Press,
  Cambridge, 1994.

\bibitem{Zwarich_MSc}
{\sc C.~Zwarich}, {\em Von {N}eumann algebras for harmonic analysis}, Master's
  thesis, University of Waterloo, 2008.

\end{thebibliography}
